\documentclass[12pt]{article}
\usepackage{euscript}
\usepackage{a4wide}
\usepackage[T2A]{fontenc}
\usepackage[utf8]{inputenc}
\usepackage{amsfonts}
\usepackage{amssymb, amsthm}
\usepackage{amsmath}
\usepackage{mathtools}
\usepackage{needspace}
\usepackage{lipsum}
\usepackage{comment}
\usepackage{cmap}
\usepackage[pdftex]{graphicx}
\usepackage{hyperref}
\usepackage{epstopdf}
\usepackage[matrix, arrow, 	curve]{xy}
\usepackage{amscd}
\usepackage{esint}
\usepackage{authblk}
\usepackage{multicol}
\usepackage{float}
\usepackage[font=small,skip=0pt]{caption}

\usepackage{wrapfig}

\begin{document}

\renewcommand{\proofname}{Proof}

\renewcommand{\d}{\partial} 
\newcommand{\Z}{\mathbb{Z}}
\newcommand{\N}{\mathbb{N}}
\newcommand{\R}{\mathbb{R}}
\newcommand{\Q}{\mathbb{Q}}
\newcommand{\K}{\mathbb{K}}
\newcommand{\Cm}{\mathbb{C}}
\newcommand{\Pm}{\mathbb{P}}
\newcommand{\B}{\mathcal{B}}
\newcommand{\Zero}{\mathbb{O}}
\newcommand{\ilim}{\int\limits}
\newcommand{\slim}{\sum\limits}
\newcommand{\action}{\curvearrowright}
\newcommand{\E}{\mathbb{E}}
\newcommand{\BB}{\overline{B}}
\newcommand{\D}{\mathcal{D}}
\newcommand{\T}{\mathcal{T}}
\newcommand{\F}{\mathcal{F}}
\newcommand{\Sf}{\mathbb{S}}
\newcommand{\Vol}{\mbox{V}}
\newcommand{\mint}{\strokedint\limits}
\newcommand{\const}{\mbox{const}}
\newcommand{\supp}{\mbox{supp}}
\newcommand{\dist}{\mbox{dist}}
\newcommand{\Hess}{\mbox{Hess}}
\newcommand{\Ker}{\mbox{Ker}}
\newcommand{\Hd}{\mathcal{H}}
\renewcommand{\div}{\mbox{div}}
\newcommand{\diam}{\mbox{diam}}
\newcommand{\NTlim}{\mbox{n.t.lim }}
\newcommand{\mydet}{\mbox{det}}
\newcommand{\Id}{\mbox{Id}}
\renewcommand{\Re}{\mbox{Re}}
\renewcommand{\Im}{\mbox{Im}}

\theoremstyle{plain}
\newtheorem{thm}{Theorem}
\newtheorem{lm}{Lemma}
\newtheorem*{st}{Statement}
\newtheorem*{prop}{Properties}
\newtheorem*{cl}{Claim}

\theoremstyle{definition}
\newtheorem*{defn}{Definition}
\newtheorem{ex}{Ex}
\newtheorem{cor}{Corollary}

\theoremstyle{remark}
\newtheorem{rem}{Rem}
\newtheorem*{note}{Note}

\title{Good elliptic operators on snowflakes}
\author[]{Polina Perstneva \thanks{The author was partially supported by the Simons Foundation grant 601941, GD.} \\ 
}
\date{\today}

\maketitle

\begin{abstract}

    We construct elliptic operators with scalar coefficients on the complements $(\R^2 \setminus S)^+$ of some Koch-type snowflakes $S$, whose Hausdorff dimensions cover the full range $(1, \ln{(4)}/\ln{(3)})$, such that the operator's elliptic measures are 
    equal to the Hausdorff measure on the boundary. This provides another example of the phenomenon that, though purely unrectifiable boundaries of domains are often characterised by the harmonic measure being singular with respect to the Hausdorff measure on the boundary, for some purely unrectifiable boundaries one can construct an elliptic operator whose elliptic measure behaves in a drastically different way. Plus, in $\R^2$, this operator can be chosen in a way that its coefficient is scalar, as opposed to a $2 \times 2$ matrix-valued one.





\end{abstract}


\tableofcontents

\section{Introduction}


The search for necessary and sufficient geometric conditions on a domain for which the harmonic measure is absolutely continuous with respect to the Hausdorff measure on the boundary started in the beginning of 20th century. They were finally articulated by 2016 in the celebrated work by J. Azzam, S. Hoffman, M. Mourgoglou, J.M. Martell, S. Mayboroda, X. Tolsa and A. Volberg \cite{AHM3TV}. In particular, it turned out that the purely unrectifiable sets are those for which harmonic measure is necessarily singular with respect to the Hausdorff measure of the boundary of the domain. Similar characterisation is true for the elliptic measure associated to more general elliptic operators which are in some sense close to Laplacian, see \cite{KP}, \cite{AGMT}, \cite{HMM}, and \cite{HMT}. 

However, even for some purely unrectifiable sets, one can find an elliptic operator (in divergence form) with scalar coefficient whose elliptic measure is absolutely continuous and even proportional to the Hausdorff measure on that set. The first example of such an anomaly was recently given by G. David and S. Mayboroda in \cite{DM}. 

Let us briefly describe their result. They work in $\R^2$ with the four corner Cantor set $K$ (so-called Garnett-Ivanov set), which is a one-dimensional Ahlfors regular set with harmonic measure (for the Laplacian) mutually singular to the Hausdorff measure $\Hd^1$ on $K$. They construct an elliptic operator in divergence form
\begin{equation}\label{op}
    L = -\div a \nabla,
\end{equation}
where $a$ is a continuous scalar function on $\R^2 \setminus K$ (as opposed to a general matrix-valued coefficient) such that 
\begin{equation}\label{scalar_coeff}
    C^{-1} \le a(x) \le C \quad \mbox{for} \; x \in \Omega = \R^2 \setminus K,
\end{equation}
and such that its elliptic measure $w_L^\infty$ with pole at infinity is equal to the Hausdorff measure $\Hd^1|_K$ (for now think of $w_L^\infty$ as sort of a universal object taken not to care about the pole). Consequently, for all poles $x$ far away from $K$, one has $C^{-1}\Hd^1|_K \le w^x_L \le C\Hd^1|_K$. For a more precise statement see Section 4 of \cite{DM}. Recall that the elliptic measure $w_L^x$ of the operator $L = -\div A \nabla$ is the measure on the boundary $\partial \Omega$ of the domain $\Omega$ which provides an integral representation of the solution of the Dirichlet problem
\begin{equation}
    \left\{\begin{matrix}
        L u = 0 \; \mbox{in} \; \Omega,\\
        u = g\ \; \mbox{on} \; \partial \Omega,
        \end{matrix}\right.
\end{equation}
namely, if $g$ is a continuous function with compact support on $\partial \Omega$, then for $x \in \Omega$ we have
\begin{equation}
    u(x) = \ilim_{\partial \Omega}{g(y)d w_L^x(y)}.
\end{equation}
For more information about the elliptic measure, see, for example, \cite{Da2}, or \cite{KKPT} and references therein. 

An important point here is that, even though the question of existence of an elliptic operator for which the corresponding elliptic measure is absolutely continuous with respect to the Hausdorff measure on the Cantor set was open for general operators in divergence form $L = -\div A \nabla$, the authors of \cite{DM} set a goal to construct an operator with a scalar-valued coefficient $a$ (in place of a matrix-valued $A$). Such a restriction seems to be more geometrically relevant, evidently makes the problem harder, and simultaneously opens the path to other interesting questions. 

Indeed, consider the similar question where the set $K$ in \eqref{scalar_coeff} is a Koch-type snowflake of dimension $d > 1$ in $\R^2$ instead of the Cantor set (for example, the usual Koch snowflake with $d = \frac{\ln(4)}{\ln(3)}$, built by an iterative procedure from an equilateral triangle). Then an elliptic operator $L = -\div A \nabla$ with a matrix-valued coefficient such that the measure $w_L^x$ is absolutely continuous with respect to the Hausdorff measure $\Hd^d|_K$ (denoted later by $w_L^x << \Hd^d|_K$) is provided immediately by a quasiconformal mapping $f: \R^2 \to \R^2$ that maps the line $l$ or a circle to $K$. One just has to use $f$ to transfer the Laplacian on a component of $\R^2\setminus l$ to a corresponding component of $\R^2 \setminus K$. The general theory states that the resulting operator $L$ is indeed an elliptic operator in divergence form, $L = -\div A \nabla$, and $w_L << \Hd^d|_K$ is true just because the same holds true for the Laplacian and the Lebesgue measure on the line. But, first, we neither have a guarantee that the coefficient $A$ is scalar, nor we know that it has any regularity at all. Second, clearly the same trick cannot be repeated with the four-corner Cantor set. This indicates that the class of operators with scalar coefficients, as opposed to matrix-valued ones, is more reasonable to study. 

This leads us to the following questions. For at least one set of some dimension $d > 1$ on $\R^2$, for example, for a Koch-type snowflake, does a good operator \eqref{op} (in the same sense as above that $w_L^x << \Hd^d|_K$) with a scalar coefficient exist? Does such an operator exist in $\R^2$ for every snowflake of dimension $d > 1$? For some other reasonable bigger class of fractals? For all decent (for example, Ahlfors regular) unrectifiable sets? 

In this paper we give a positive answer for the first of these three questions for all $d$ such that $1 < d < \frac{\ln(4)}{\ln(3)}$. The last question is a very ambitious one, and we have no idea if a satisfactory solution to that problem can exist, but, concerning the second question, we believe that one can build the desired operators for a considerably larger class of fractals than described in this paper. For perspectives of generalisations and extensions see Section 8.

Our result is also connected to the questions about the dimension of harmonic (or elliptic) measure and about the dimension drop. For the harmonic measure in $\R^2$, by the result of Wolff \cite{W1} dated 1993, preceded by works of Makarov and Jones and Wolff, the dimension of the harmonic measure is necessarily one. For $\R^n$ with $n \ge 3$ the question appears to be much more complicated. On the one hand, Wolff constructed examples of domains $\Omega_n \subset \R^n$ such that the dimension of harmonic measure $\dim w_{\Omega_n}$ is strictly larger than $n - 1$, \cite{W2}. On the other hand, the dimension of harmonic measure can never reach $n$, due to Bourgain, \cite{B}: for every $n \ge 3$ there exists a universal constant $b_n > 0$ such that for any domain $\Omega$ we have $\dim w_\Omega \le n - b_n$. Bourgain's constant $b_n$ is not optimal, and is very hard to compute explicitly, and one of the most celebrated and challenging open problems connected to harmonic measure up-to-date is to find the optimal value of $b_n$. The only significant progress made in the last 30 years there, up to the authour's knowledge, is a recent work by Badger and Genschaw, who provided an explicit lower bound on $b_n$ by revising Bourgain's proof, \cite{BG}. Bishop's conjecture that the optimal value of $b_n$ is $\frac{1}{n - 1}$, still stands, \cite{Bi}. For the elliptic measure, without any regularity assumptions on the coefficient of the corresponding elliptic operator, things can differ even in dimension $2$ already. Thus, for every $\varepsilon > 0$ one can construct a planar domain and an elliptic operator $L$ in divergence form such that the associated elliptic measure has dimension at least $2 - \varepsilon$, see Sweezy \cite{S}. 

One says that a dimension drop occurs when the dimension of the harmonic (elliptic) measure $\dim w_\Omega$, which satisfies naturally $\dim w_\Omega \le \dim \partial \Omega$, is strictly less than the dimension of the boundary of the domain $\Omega$. For the harmonic measure, there are plenty of examples when this happens, including some Wolff snowflakes, some self-similar sets, and even sets with no nice geometric structure (in the classical sense): see a work by Azzam \cite{Az} and references wherein. The analogue of the Azzam's result is also likely to be true for operators close to the Laplacian. 

Our result implies that for elliptic measures on Koch-type snowflakes associated to our operators a dimension drop does not occur. Also, though the absolute continuity of the elliptic measure with respect to the $d$-dimensional Hausdorff measure $\Hd^d$ on the boundary is (much) stronger than the property that the elliptic measure has dimension $d$, our result can be interpreted as a stronger version of the result by Sweezy for operators with continuous scalar-valued coefficients, as opposed to general matrix-valued elliptic coefficients, in the range of dimensions $(1, \ln{(4)}/\ln{(3)})$. It would be, once again, interesting to know whether it is possible to cover the whole range $(1, 2)$.

Before we state our main theorem, we describe the sets $K$ (the family of snowflakes) we will work with. Because of personal preferences, our snowflakes are going to be non-compact (though we can do the compact versions as well, see Section 8). Also, we use a slightly different replacement algorithm than the one for the classical Koch snowflake. In order to not confuse the reader, we'll keep the precise definition with formulas for later, but we sketch the construction now. We start from the unit interval $I = [-1/2, 1/2]$. Fix any $\alpha, 0 < \alpha < \pi/3$. Replace the segment $I$ with $F(I) = F_\alpha(I)$ which is constructed as follows. (Actually, we will give later a more qualified definition of the transformation $F$, which will play the crucial role in our construction.) Replace the middle of the interval by two sides of the isosceles triangle built above $I$, with the angle $\alpha$ between the leg and the base, so that the four intervals we obtain after this procedure have equal length, and that the picture stays symmetric with respect to the bisection of $I$:

\begin{figure}[H]
    \centering
    \includegraphics[scale = 0.2]{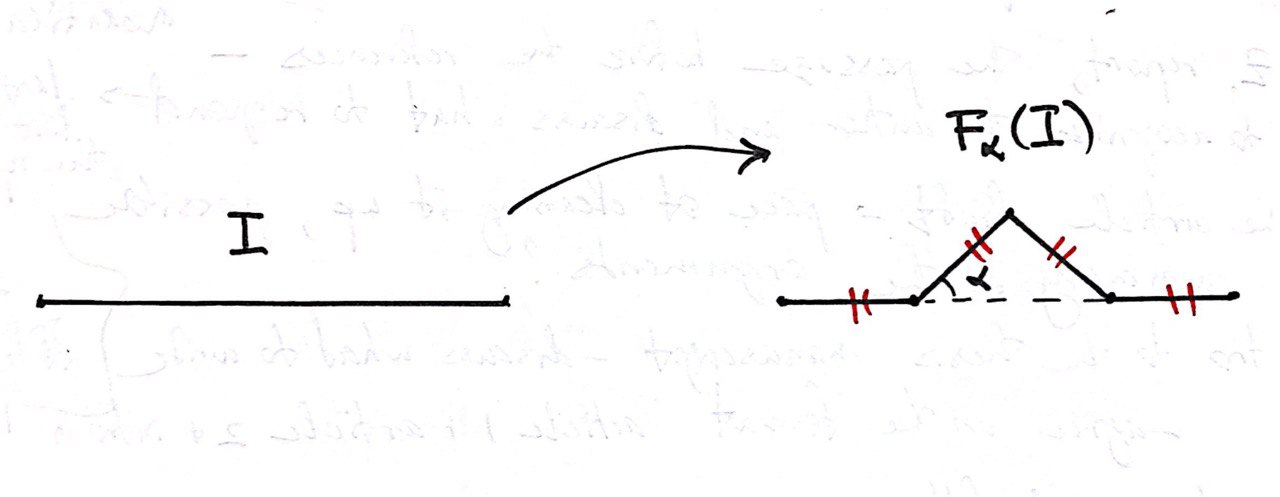}
    \caption{Simple version of transformation F}
    \label{fig:F alpha I v one}
\end{figure}

Repeat the procedure for each of the four intervals composing $F(I)$ in place of $I$, and iterate. Concerning how our sets $K$ look globally, we consider two examples of Koch-type snowflakes $S_\alpha^1$ and $S_\alpha^2$ which are locally the same, but $S^1_\alpha$ is purely fractal, and $S^2_\alpha$, on a large scale, resembles a line. To get $S_\alpha^2$, we decompose $\R = \cup_{k \in \Z}{[k - 1/2, k + 1/2]}$, and for all $k \neq 0$ repeat the procedure we applied to the unit interval $[-1/2, 1/2]$. 

\begin{figure}[H]
    \centering
    \includegraphics[scale = 0.3]{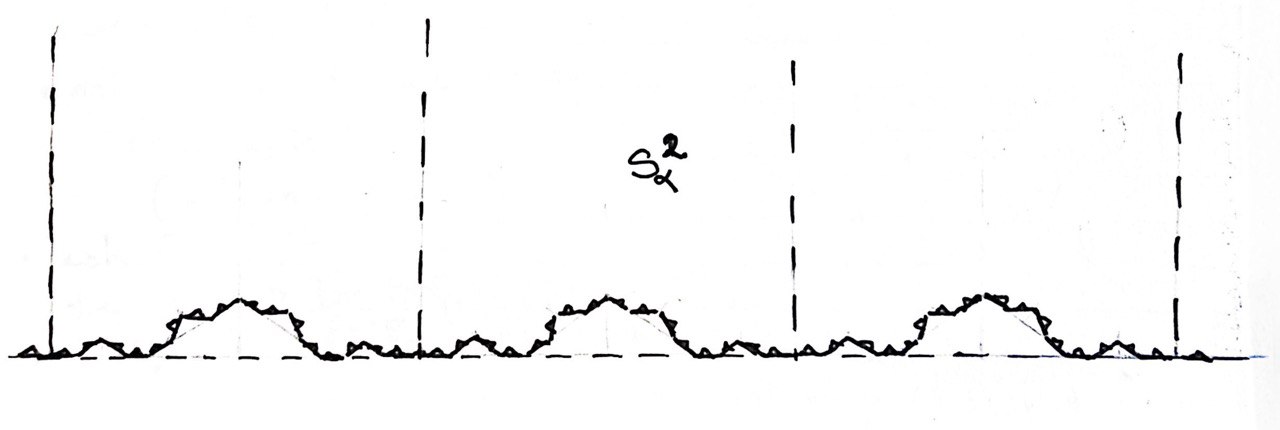}
    \caption{$S^2_\alpha$}
    \label{fig:S 2 alpha}
\end{figure}

To get $S_\alpha^1$, first denote by $l$ the length of any segment the transformation $F$ described above produces from the unit interval $I$. It is easy to compute that $l = \frac{1}{2(1 + \cos{\alpha})}$. By $H_l$ denote the scaling transform (homothety) with the coefficient $1/l$ and center $-1/2$ (the left end of the unit interval $I$), and denote by $S_\alpha |_I$ the fractal set we built from $I = [-1/2, 1/2]$. Define the part of the set $S^1_\alpha$ to the right of $\{x = -1/2\}$ to be $\cup_{k \ge 1}{H_l(S_\alpha |_I)}$. To get the part of the set to the left of the axis, we reflect the part to the right with respect to it: $S^1_\alpha = \cup_{k \ge 1}{H_l(S_\alpha |_I)} \cup \left(\overline{\left(\cup_{k \ge 1}{H_l(S_\alpha |_I)} + 1/2\right)} - 1/2\right)$. Here is how $S^1_\alpha$ looks like to the right of $-1/2$.


\begin{figure}[H]
    \centering
    \includegraphics[scale = 0.25]{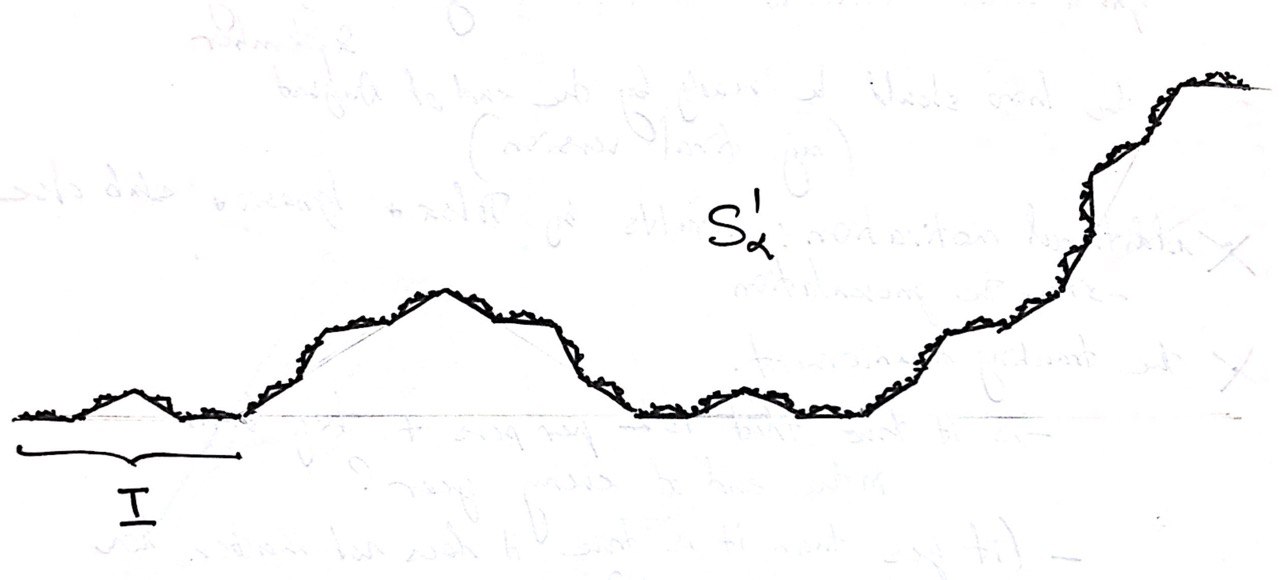}
    \caption{$S^1_\alpha$}
    \label{fig:S 1 alpha}
\end{figure}

In the theorem below we will call $S_\alpha$ both $S_\alpha^1$ and $S_\alpha^2$ for simplicity. A classical computation, essentially the same as for the standard Koch snowflake, shows that the Hausdorff dimension $d$ of $S_\alpha$ (both $S_\alpha^1$ and $S_\alpha^2$) is $\frac{\ln{(4)}}{\ln{(2(1 + \cos{\alpha}))}}$. We also choose the domains $\Omega_1 = \left(\R^2 \setminus S_\alpha^1\right)^{+}$ and $\Omega_2 = \left(\R^2 \setminus S_\alpha^2\right)^{+}$, and denote them both $\Omega$ in what follows. 

We now state our main result.

\begin{thm}\label{mainthm1}
    For each symmetric Koch-type snowflake $S_\alpha$ with $0 < \alpha < \pi/3$ there exists an elliptic operator $L = -\div \, a \nabla $ with a continuous scalar coefficient $a$ such that 
    \begin{enumerate}
        \item $a$ is defined on the domain $\Omega$ and is continuous,
        \item as in \eqref{scalar_coeff}, $a$ is bounded and bounded away from zero,
        \item and the elliptic measure $w^\infty_L$ (to be defined below) is equal to the Hausdorff measure $\Hd^d|_{S_\alpha}$, and, in addition, $C^{-1}\Hd^d|_{S_\alpha} \le w^x_L \le C\Hd^d|_{S_\alpha}$ for all $x$ such that $\delta(x) = \dist (x, S_\alpha) \ge 1$, where $1 < d = \frac{\ln{(4)}}{\ln{(2(1 + \cos{\alpha}))}} < \frac{\ln(4)}{\ln(3)}$.
    \end{enumerate} 
\end{thm}

\begin{rem}
    For each snowflake $S_\alpha$ we build a different elliptic operator, or, a different coefficient $a$, so the correct notation would be $L_\alpha = -\div \, a_\alpha \nabla$ and $w^x_{L_\alpha}$. However we always work with one snowflake, or a fixed $\alpha$, so we usually omit this dependence in the notation. The constant $C$ in $\eqref{scalar_coeff}$ and in $C^{-1}\Hd^d|_{S_\alpha} \le w^x_L \le C\Hd^d|_{S_\alpha}$ in property $3$ above also depends on $\alpha$, and it will be clear from the construction that we get a constant $C = C(\alpha)$ that tends to infinity as $\alpha$ tends to $\pi/3$. This is connected to why the case of the classical Koch snowflake ($\alpha = \pi/3$) is not covered by Theorem $\ref{mainthm1}$. See a short discussion in Section 8 about this.
\end{rem}

\begin{rem}
    Regarding the regularity of our scalar coefficient $a$, we cannot guarantee anything better than continuous on $\Omega$. We cannot get H\"older continuous $a$, since, as we get closer to the boundary, $a$ oscillates more and more, which will be clear from the construction. This accords with the perturbation results about the behaviour of elliptic measures of operators ``close to'' Laplacian. Indeed, locally H\"older continuous coefficients satisfy, for example, the condition $(1.3)$ in \cite{AGMT}. So, the behaviour of elliptic measure of an operator with locally H\"older continuous coefficient, in what concerns the relations between rectifiability of the boundary and absolute continuity with respect to the boundary measure, resembles the behaviour of the harmonic measure. In particular, the elliptic measure of an operator with locally H\"older continuous coefficient is singular with respect to $\Hd^d|_{S_\alpha}$. Therefore we get again that our coefficient $a$ cannot be locally H\"older continuous. 

    Speaking of different regularity conditions on the coefficients, and what happens in between the ``standard'' case of operators ``close'' to Laplacian and the dramatically different case of operators like we construct, some results in this direction for the planar case were recently obtained by I. Guillén-Mola, M. Prats, and X. Tolsa, \cite{GPT}. In particular, they proved the following: if the coefficient $A$ of an elliptic operator in the divergence form is Lipschitz, and the boundary of a planar domain $\partial \Omega$ is a Reifenberg flat set with a small constant (how small it should be depends on the ellipticity of $A$), then there is a subset of the boundary of $\sigma$-finite length such that it has elliptic measure one. They also classify settings by the regularity of the coefficient $A$ and the domain, giving in each case a lower estimate on the dimension of the elliptic measure.
\end{rem}

\medskip

In the core of the construction of our operators lies a restoration procedure of the scalar coefficient $a$ from the level lines of a pair of (conjugated) functions $(u, v)$ which solve the equations $-\div \, a \nabla$ and $-\div \, a^{-1} \nabla$ respectively, see Section 2 both here and in \cite{DM}. Without entering into the details here, the coefficient $a$ is reconstructed from a net of two families of mutually orthogonal curves (any curve from one family is orthogonal to any curve from the second one) with a special property that ``each cell of the net does not differ too much from a square uniformly on $\Omega$''. Thus, our main goal is to build such a net of curves. Our net construction resides on a tiling idea. First, we will build a family of partially self-similar curves $s_k$ which approximate our snowflake $S_\alpha$. Then we use these curves to divide the domain $\Omega$ above $S_\alpha$ into stripes between two consecutive $s_k$ and chop each stripe into symmetric tiles of two types, getting a large scale net. We fill insides of our tiles with local nets with the property we need, using the symmetry of tiles. Finally we check that local nets glue nicely, and that the special net property for the whole $\Omega$ follows, because we use a compact set of tiles to cover the whole space, and build the large scale net as to, vaguely, satisfy the same special property. 

In Section 3 we describe in details our net construction for $S_\alpha$. In Section 4 we explain why the cells for this net tile the space. In Section 5 we smoothe out our net a bit, since the smoothness is required by the reconstruction procedure described in the previous paragraph. In Section 6 we finally draw the intermediate curves and check the special property mentioned in the previous paragraph. In Section 7 we state our main theorem more precisely and finally finish its proof. In Section 8 we make some remarks, state some hypotheses and explain why they might be (hopefully) interesting and why they are difficult to tackle for us up-to-date.

\subsubsection*{Acknowledgements} 

I am very grateful to my PhD thesis advisor, Professor Guy David, for introducing me to the problem, and for fruitful discussions. Many thanks both to him and to Antoine Julia for reading and making useful suggestions to improve the earlier versions of this manuscript.

\section{Building the coefficient $a$ from the level lines of a Green function}

In this section we describe how to construct (locally) the scalar coefficient $a$ from a function $u$ in such a way that $\div a \nabla u = 0$. The procedure is somewhat geometric, and will rely mostly on the level lines of the function $u$. Furthermore, it will be explained how certain properties of the family of level lines of $u$ are translated to the properties of the coefficient $a$. Therefore, in a general perspective, our goal will be to construct a function $u$, which will end up being the Green function with the pole at infinity (a universal solution in some sense of the equation $\div a \nabla \cdot = 0$), with certain desired properties. So that, at the end of the day, the restored coefficient $a$ will satisfy all the conditions in Theorem \ref{mainthm1}. 

This part of the proof of Theorem \ref{mainthm1} is the same as for the analogous result of David and Mayboroda about the four-corner Cantor set, and this Section repeats Section 2 in \cite{DM} for the sake of completeness and probably makes some parts of it a little clearer. 

Let us work, to start with, in an open neighbourhood $U$ of a fixed point $y$, where a function $u$ such that $u \in C^4$ and $\nabla u \neq 0$ is defined. Consider the set of level lines of the function $u$, $\gamma_s = u^{-1}(s)$: by the implicit function theorem, each $\gamma_s$ is a piece of a smooth curve with no self-intersections (if that is not true we can always restrict ourselves to a smaller neighbourhood of the point $y$). The gradient vector field $\nabla u$ is, as it is well-know, orthogonal to the level lines of $u$: for every $x \in U$, if $x \in \gamma_s$, then $\nabla u(x)$ is orthogonal to the tangent vector to $\gamma_s$ at $x$. Using $\nabla u$, we build a family of curves orthogonal to $\gamma_s$ in the following way. Without loss of generality, the image $u(U)$ contains $0$. We parameterise $\gamma_0$, $\theta \to \gamma_0(\theta)$ (locally, since we still work in $U$), and solve locally the equation $z'(t) = \nabla u(z(t))$ with the initial condition $z(0) = \gamma_0(\theta)$. The solution of this equation is a $C^3$ curve, which can be extended in both direction of $t$ as long as it stays in $U$. We will call this maximal extension $\Gamma_\theta$. The standard theorem about the existence and uniqueness of solutions to differential equations of the first order asserts that $\Gamma_{\theta_1}$ and $\Gamma_{\theta_2}$ do not cross as long as $\theta_1 \neq \theta_2$.

By definition the $\gamma_s$ cover $U$. We claim that the curves $\Gamma_\theta$ also cover $U$ if it is small enough: indeed, at any point of any $\gamma_s$ there exists an orthogonal curve $\Gamma_\theta$ which passes through it. Again, without loss of generality, assume that $y \in \gamma_0$ and that $y \in \Gamma_{\theta_0}$. Then, if $z = z(t) \in \Gamma_{\theta_0}$ and close enough to $y$, every point $z'$ which is close enough to $z$ lies on some $\Gamma_\gamma$. Indeed, let us solve the equation $w(r) = -\nabla u (w(r))$ with the initial data $w(0) = z'$. Since $z'$ is close enough to $z$, theorem about good dependence on the initial data and parameters of the solution of a first-order differential equation guarantees that the curve $w$ crosses $\gamma_0$ at a point $\gamma_0(\theta)$, close to $y$, at a moment $r \approx t$. This is true because the gradient field $-\nabla u$ along $w$ stays close to the gradient field $-\nabla u$ along the curve $\Gamma_{\theta_0}$ ran backwards, and that the latter curve $\Gamma_{\theta_0}$ is the solution for the same equation $w(r) = -\nabla u (w(r))$ with the initial data $w(0) = z$ (which crosses $\gamma_0$ at $y$). So, we can define a mapping $v$ on $U$ (near $y$) which to $z'$ associates the unique $\theta$ such that $z'$ lies on $\Gamma_\theta$. This mapping $v$ is at least $C^2$, and $\nabla v \neq 0$ near $y$ thanks to the implicit function theorem. (Here we mean that there is a formula which connects $z'$ with $\theta$: we know that there exists $r$ such that $z' = z'(r) = \gamma_0(\theta) + \ilim_0^r{\nabla u(z'(s))ds}$, and the formal derivative of the expression with respect to $\theta$ is nonzero. Therefore $\theta$ is a function of $z'$, $\theta = v(z')$, with nonzero derivative $\nabla v \neq 0$.) 

Thus we just reconstructed from $u$ a pair of functions $(u, v)$ which can be used as coordinates near the point $y$ (in the open neighbourhood $U$). By definition, the pair $(u, v)$ satisfies the orthogonality relation
\begin{equation}\label{perp}
    \nabla u(z) \perp \nabla v(z).
\end{equation}
We remark here that the change of variables defined by $(u, v)$ is not necessarily conformal, because for this vectors $\nabla u$ and $\nabla v$ would have to have the same length, and this would be the case only if the initial function $u$ was harmonic.

Now we show that in any open set where $u$ and $v$ are well-defined, of class $C^2$, such that $\nabla u, \nabla v \neq 0$ and \eqref{perp} holds, we have
\begin{equation}\label{the_equation}
    \div \, a \nabla u = 0 \quad \mbox{with} \; a(z) = \frac{|\nabla v(z)|}{|\nabla u(z)|},
\end{equation}
and $\div \, \frac{1}{a} \nabla v = 0$ with the same $a$. 
Indeed, observe that the vector $w = \left(\frac{\partial v}{\partial y}, -\frac{\partial v}{\partial x}\right)$ is orthogonal to $\nabla v$ and has the same length as $\nabla v$. At the same time, $w$ is proportional to $\nabla u$ because of \eqref{perp} and the fact that we work in dimension 2. Thus, $a \nabla u = \pm w$ with the coefficient $a$ defined as above. Since in $a \nabla u = \pm w$ the sign on the right-hand side is locally constant, we have
$$\div \, a \nabla u = \div \, w = \frac{\partial^2 v}{\partial x \partial y} - \frac{\partial^2 v}{\partial x \partial y} = 0,$$
as in \eqref{the_equation}. 

Up to now we only worked locally. But of course, if our function $u$ of class $C^4$ is defined and satisfies $\nabla u \neq 0$ everywhere on a domain $\Omega$, or the described above coordinate lines are already given to us on the whole domain, the above provides us with the desired coefficient $a$ everywhere on $\Omega$. This will be our case.

Moreover, in this article we will be able to reconstruct the values of the parameters $\theta$ and $s$ from a level curve of $v$ or $u$ respectively. We stress that later we will rely on that fact. 

\begin{rem}\label{rem1}
We have to point out on the side however that, in general, on the global picture level curves are not curves in the strict sense of the word: level sets can have several connected components (which are curves), which can be cyclic. See, e.g., the main example in \cite{DM}, p. 15 or 17. Therefore, one cannot always reconstruct precisely a level curve from a value $s$ (of a function $u$) or $\theta$ (of a function $v$). This does not contradict that the coefficient $a$ is well-defined globally. It just means that, to enlist all the (important for the construction) level curves of $v$ by parametrizing one level curve $\gamma_0$ of $u$, one should get very lucky with the choice of $\gamma_0$ (which was the case in \cite{DM}). We will be lucky, and even better: for all values of $s$ and $\theta$ we will have the unique level curve of $u$ or $v$ respectively corresponding to it (and each curve will have infinite length). A choice of the level curve $\theta_0$ will be made in advance, but afterwards it will be clear that we could have chosen any level curve of $u$. 
\end{rem}

~\

Thus we have built the scalar coefficient $a$ from a function $u$ so that the latter is a solution to $\div \, a \nabla \cdot = 0$ on $\Omega$. Recall however that our interest is also to ensure the double inequality $C^{-1} \le a \le C$. The rest of the section will be devoted to reformulating and specifying what needs to be checked in order to do this, and we will use the results of this preparatory work starting from Section 6.

The double inequality $C^{-1} \le a \le C$ by $\eqref{the_equation}$ is the same as the double inequality 
\begin{equation}\label{double_ineq}
    C^{-1}|\nabla u| \le |\nabla v| \le C|\nabla u|.
\end{equation}
Our aim throughout the upcoming proof of Theorem $\ref{mainthm1}$ is to control that $\eqref{double_ineq}$ is fulfilled, while we construct the level lines of $u$ and $v$. In terms of those lines, assuming we fixed a nice parameterization $\theta \mapsto \gamma_0(\theta)$, the gradient of $u$ is locally proportional to an increment $\delta_s$ divided by the distance between the level sets. Similarly, the gradient of $v$ is proportional to the inverse of the distance between the $\Gamma_\theta$ curves, divided by the increment $\delta_\theta$. To control $a$, we ensure that the ratio of these two quantities stays between constants. This means that, on the picture, if we look at a cell of our coordinate system obtained from roughly equal increments $\delta_s$ and $\delta_\theta$, the cell should resemble a square, or a rectangle which is not too thin. Notice that nothing bad happens if the distance between the level lines of $u$ varies a lot, as long as in the same place the same happens with the level lines of $v$. 

More precisely, assume that we want to estimate the coefficient $a$ at a point $z \in \R^2 \setminus K$. In what follows we will call level lines of $u$ green lines, and level lines of $v$ -- red lines, so as to be able to distinguish between the two more easily. Consider a ``cell'' containing $z$ bounded by green lines, which correspond to values $s_1$ and $s_2$ of the function $u$, and two red lines, which correspond to values $\theta_1$ and $\theta_2$. Here we just postulate the correspondences between the curves and the values, but, for later, see Remark \ref{rem1} about this reconstruction (or, we can also think that the cell considered is small enough, and locally the reconstruction considered is always possible). 

\begin{figure}[H]
    \centering
    \includegraphics[scale = 0.08]{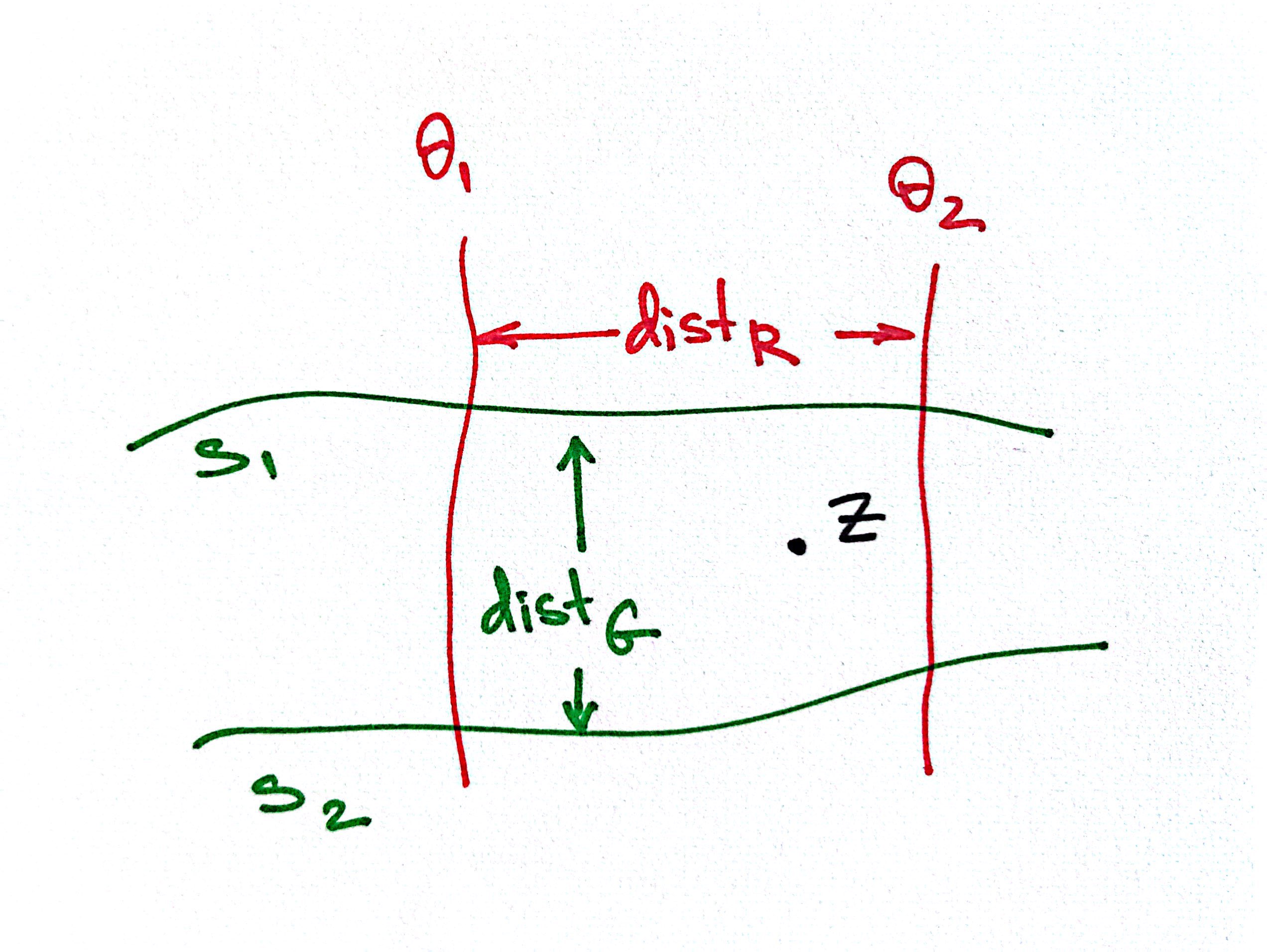}
    \caption{A cell of the net of green and red curves}
    \label{fig:cell}
\end{figure}

The cells will be regular enough; the point will be to control the ratio of distances between the sides, as follows. Let $\dist_G$ be the distance between the green sides of the cell and $\dist_R$ be the distance between the red sides of the cell. By ``green sides of the cell'' we mean the two pieces of green curves $s_1$ and $s_2$ which are bounded by the red curves $\theta_1$ and $\theta_2$, so $\dist_G$ is the distance between these pieces, not between the curves $s_1$ and $s_2$. Similarily, $\dist_R$ is the distance between two pieces of $\theta_1$ and $\theta_2$, bounded by the green curves $s_1$ and $s_2$. If the size of the cell ($\dist_G$ and $\dist_R$) is small enough, then 
$$a(z) = \frac{|\nabla v|}{|\nabla u|} \le c_0 \frac{1/\dist_R \delta_s}{1/\dist_G \delta_\theta} = c_0\frac{|s_1 - s_2|\dist_G}{|\theta_1 - \theta_2|\dist_R}, \quad \mbox{and similarily}$$
\begin{equation}
   c_0^{-1}\frac{|s_1 - s_2|\dist_G}{|\theta_1 - \theta_2|\dist_R} \le a(z),
\end{equation}
where $c_0$ is a constant close to $1$ depending on the size of the cell, which gets closer to $1$ as the cell gets smaller. Thus, if we will ensure the double inequality
\begin{equation}\label{mainproperty}
    C^{-1} \le \frac{|s_1 - s_2|\dist_G}{|\theta_1 - \theta_2|\dist_R} \le C 
\end{equation}
for all $z \in \R^2 \setminus K$ for all small enough (depending on $z$) cells (small $|s_1 - s_2|$ and $|\theta_1 - \theta_2|$), we will ensure the desired property $\eqref{scalar_coeff}$ of the coefficient $a$.

\begin{rem}
    We mentioned that we are going to build the net of level lines of $u$ and $v$. Of course, it is enough to know the level lines of $u$, because from them we can deduce the direction of the gradient and draw the level lines of $v$. But we will draw both families of lines more of less at the same time, because it will help us to guarantee \eqref{scalar_coeff} at each step of our construction. 

    The picture stays the same if we relabel the level lines of $u$, that is, replace the function $u$ with a function $f \circ u$ (for a reasonable $f$), or if we change the parameterisation of the ``zero'' curve $\gamma_0$, that is, replace $v$ by $g \circ v$.

    This evidently changes the coefficient $a$ a bit, but, if we control the size of $f'$ and $g'$, it does not spoil the property we want from it. When $u$ satisfies $\div \, a \nabla = 0$, a function $w = f \circ u$ satisfies the equation $\div \, b \nabla w = 0$ with $b(z) = a(z)/f'(u(z))$, since $b(z)\nabla(f \circ u)(z) = b(z) f'(u(z)) \nabla u(z) = a(z) u(z)$. If $f'$ is bounded and bounded away from zero (which characterises a reasonable relabelling or a reasonable reparameterisation), then \eqref{scalar_coeff} is still satisfied with another constant $C$. The same happens if we change $v$, the conjugated function which satisfies $\div \, \frac{1}{a}\nabla$, for $g \circ v$. One can say that a relabeling or a reparameterisation, which keeps the picture the same, leaves the coefficient $a$ in the same class of equivalence we are studying.     

\end{rem}

\section{Building the large scale net for a compact subset of the snowflake}

In this section we start building the net of curves to reconstruct the scalar coefficient $a$, as explained in the previous section. For now, we do this in a neighbourhood of the part of the fractal generated by the unit interval $[-1/2, 1/2]$, as in figure $\ref{fig:F alpha I v one}$. It will be the main step of our construction: covering the rest of the domain $\Omega$ later will be a question of iterating the transformation $F$ we'll define in a minute. First, we give the accurate inductive definition of the part of the snowflake $S_\alpha|_{[-1/2, 1/2]}$, which interests us now, and construct a discrete family of green curves labelled $\{s_k\}$ (here we use again that we can identify a green curve and a value of the function $u$ on it, and we label curves with these values; the exact numbers $s_k$ are to be precised later). These curves will approximate the snowflake $S_\alpha|_{[-1/2, 1/2]}$ ``from above'' better as $k$ goes to infinity (in terms of the Hausdorff distance, say). Simultaneously, we construct a discrete family of red curves $\{\theta_k\}$, which together with the family $\{s_k\}$ will serve as a supporting structure (the large scale net from Introduction) to build the rest of the net (we explain later how).

Note that the two (discrete) families of green and red curves we get in this section are not going to be smooth at a discrete family of points, while in the previous section we had our curves at least $C^2$ and orthogonal. This is not critical and will be fixed later in Section 5: we will apply a mollifying procedure around those non-smoothness points, which allows us to eliminate both inconsistencies.

We now start the construction by defining the transformation $F$ we mentioned earlier. We already saw in the introduction a version of this transformation for intervals, see figure \ref{fig:F alpha I v one}, but we will need its more general version to construct $\{s_k\}$.
\begin{defn}[transformation F]
Given a parameter $\alpha$, $0 < \alpha < \pi/3$, define 
\begin{equation}\label{l}
    l = \frac{1}{2(1 + \cos{\alpha})}.
\end{equation}
Given a compact set $I$ in $\R^2 = \Cm$, points $p_I, q_I \in I$ which we will call its left end and right end, and a point $p$, define the transformation $F_\alpha(I, p)$ as
$$F_\alpha(I, p) = F(I, p) = F_1(I) \cup F_2(I) \cup F_3(I) \cup F_4(I) \quad \mbox{with}$$
$$F_1(I) = (I - p_I)l + p, \quad z_{F_1(I)} = (z_I - p_I)l + p, $$
$$F_2(I) = (I - p_I)le^{i\alpha} + q_{F_1(I)}, \quad z_{F_2(I)} = (z_I - p_I)le^{i\alpha} + z_{F_1(I)}, $$
$$F_3(I) = (I - p_I)le^{-i\alpha} + q_{F_2(I)}, \quad z_{F_3(I)} = (z_I - p_I)le^{-i\alpha} + z_{F_2(I)},$$
\begin{equation}\label{F}
     \quad F_4(I) = (I - c_I)l + q_{F_3(I)} \quad \mbox{and} \quad z_{F_4(I)} = (z_I - p_I)l + z_{F_3(I)},
\end{equation}
where $z_I$ denote any point on $I$, and $z_{F_i(I)}$ -- its images: for example, $q_{F_1(I)} = (q_I - p_I)l + p$. So the left end and the right end of the set $F(I, p)$ are also defined: $p_{F(I, p)} = p$ and $q_{F(I, p)} = q_{F_4(I)}$.
\end{defn}

\begin{figure}[H]
    \centering
    \includegraphics[scale = 0.2]{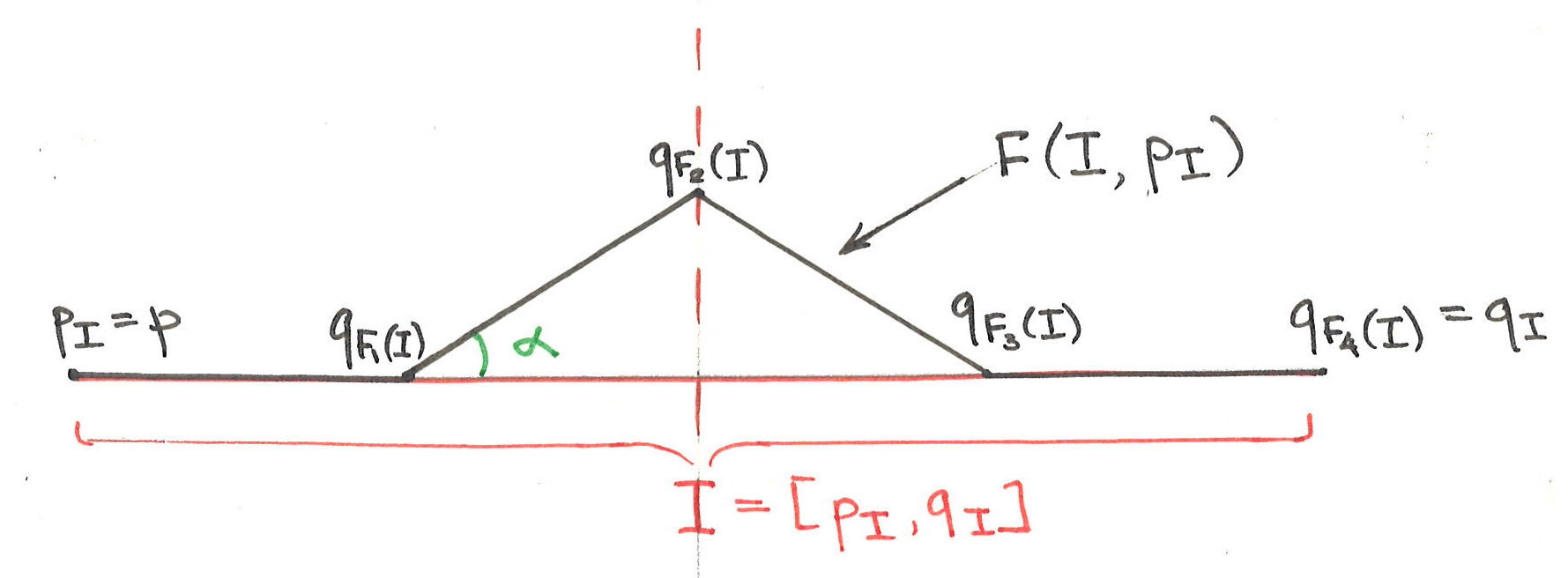}
    \caption{The advanced version of the transformation $F$}
    \label{fig:F alpha I v two}
\end{figure}

The basic example of how the transformation $F$ acts is given by $I$ being a horizontal interval $I = [p_I, q_I]$ and $p = p_I$. This is the version of the transformation mentioned in the introduction. Recall that it replaces the middle of interval $I$ with two legs of the isosceles triangle with the angle $\alpha$ and the middle of $I$ as a base. So the four intervals obtained, $[p, q_{F_1(I)}], [q_{F_1(I)}, q_{F_2(I)}], [q_{F_2(I)}, q_{F_3(I)}]$ and $[q_{F_3(I)}, q_{F_4(I)}]$, have equal length $|I|l$ and the picture stays symmetric with respect to the bisection of $I$. See figure $\ref{fig:F alpha I v two}$.

The usual inductive definition of $S_\alpha$ uses the unit interval $I = [-1/2, 1/2]$ as a base step and exploits iterations of the version of $F$ from the example applied to all the intervals obtained at the previous stage:
\begin{equation}\label{snowdef1}
    S^0_\alpha = I = [-1/2, 1/2], \quad S^k_\alpha = \cup_{I \in S^{k - 1}_\alpha}{F(I, p_I)}, \quad S_\alpha = \lim_{k \to \infty} S^k_\alpha .
\end{equation}
Given this, it should not be too surprising for the reader, especially if one is familiar enough with the Koch's snowflake, that, as often with fractals, we need to index points, intervals and arcs of curves with multi-indices. So we index ends of intervals which compose $S^k_\alpha$ and $s_k$ with words/strings of symbols $\{1, 2, 3, 4\}$ of length $k + 1$, because they encode naturally how $S^k_\alpha$ and $s_k$ are constructed. We denote these words $\mathcal{I}_k$, $\mathcal{I}_k = \{w = w_1\dots w_k w_{k + 1}, \; w_i = 1, 2, 3, 4\}$, 

Indeed, according to the definition $\eqref{snowdef1}$, from an interval $[z_{w_1}, z_{w_2}]$ at the $k$th iteration $S^k_\alpha$ (where $z_{w_1}$ and $z_{w_2}$ are its left end and right end respectively, as usual), we make four of equal length with the aid of three supplementary points $z_{w_1 2} = q_{F_1([z_{w_1}, z_{w_2}])}$, $z_{w_1 3} = q_{F_2([z_{w_1}, z_{w_2}])}$ and $z_{w_1 4} = q_{F_3([z_{w_1}, z_{w_2}])}$, where $w_1i$ is the string concatenation. We also rename $w_1$ to $w_1 1$ and $w_2$ to $w_2 1$. This allows us to reconstruct easily the location of a point $z_w$ by its index $w$: we read the word symbol by symbol and choose one of the four intervals ($w_i$th for the symbol number $i$) at each iteration, which is where $z_w$ is situated. Sometimes we say that $z_w \in S^k_\alpha$ with the word $w$ shorter than $k + 1$-symbol word. In that case we mean that $z_{w1\dots 1} \in S^k_\alpha$, where we put enough ones at the end so that $w1\dots 1$ has $k + 1$ symbols. This leads to no confusion because, by the notation above, points $z_w$ and $z_{w1\dots 1}$ are factually the same. 

Another notation remark: while the green curves $\{s_k\}$ will be enumerated by integers, in this section we enumerate the red curves $\{\theta_k\}$ with words as well for everybody's convenience, so the correct notation is $\{\theta_w\}$.

Before we finally start the construction of our green and red curves, we'd like to highlight the important role transformation $F$ plays in the whole business. While the curves $s_1$, $s_2$ and a couple of first red curves have to be defined manually, $F$ will allow us to define everything else quasi-automatically and trace the self-similarity of the construction.  

\medskip

\begin{figure}[H]
    \centering
    \includegraphics[scale = 0.15]{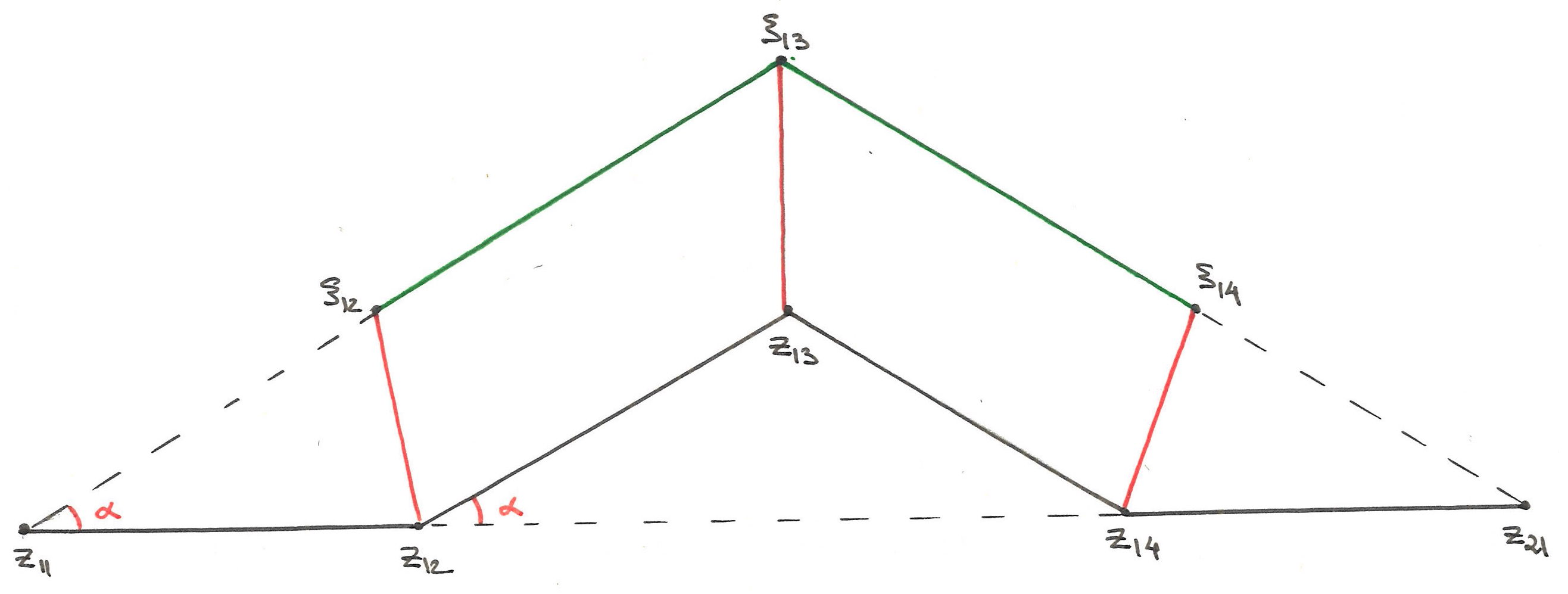}
    \caption{$S^1_\alpha$ and the curve $s_1$}
    \label{fig:snow 0}
\end{figure}

{\bfseries Step 1.} 
Recall that by the definition of $F$ we have
$$S^1_\alpha = [z_{11}, z_{12}] \cup [z_{12}, z_{13}] \cup [z_{13}, z_{14}] \cup [z_{14}, z_{21}] \quad \mbox{with}$$
$$z_{11} = z_{1}, \quad z_{12} = z_{11} + l(z_{2} - z_{1}), \quad z_{13} = z_{12} + l(z_{2} - z_{1})e^{i\alpha},$$
$$z_{14} = z_{13} + l(z_{2} - z_{1})e^{-i\alpha} = z_{21} - l(z_{2} - z_{1}), \quad z_{21} = z_{2}.$$
We define
$$s_1 = [\zeta_{12}, \zeta_{13}] \cup [\zeta_{13}, \zeta_{14}] \quad \mbox{with}$$
$$\zeta_{12} = z_{11} + l(z_{2} - z_{1})e^{i\alpha}, \quad \zeta_{13} = \zeta_{12} + \frac{l}{\cos{\alpha}}(z_{2} - z_{1})e^{i\alpha},$$
\begin{equation}\label{s1}
    \zeta_{14} = \zeta_{13} + \frac{l}{\cos{\alpha}}(z_{2} - z_{1})e^{-i\alpha} = z_{21} - l(z_{2} - z_{1}),
\end{equation}
and
\begin{equation}
    \theta_{12} = [z_{12}, \zeta_{12}], \theta_{13} = [z_{13}, \zeta_{13}], \theta_{14} = [z_{14}, \zeta_{14}].
\end{equation}
See figure $\ref{fig:snow 0}$ below.

\begin{figure}[h]
    \centering
    \includegraphics[scale = 0.08]{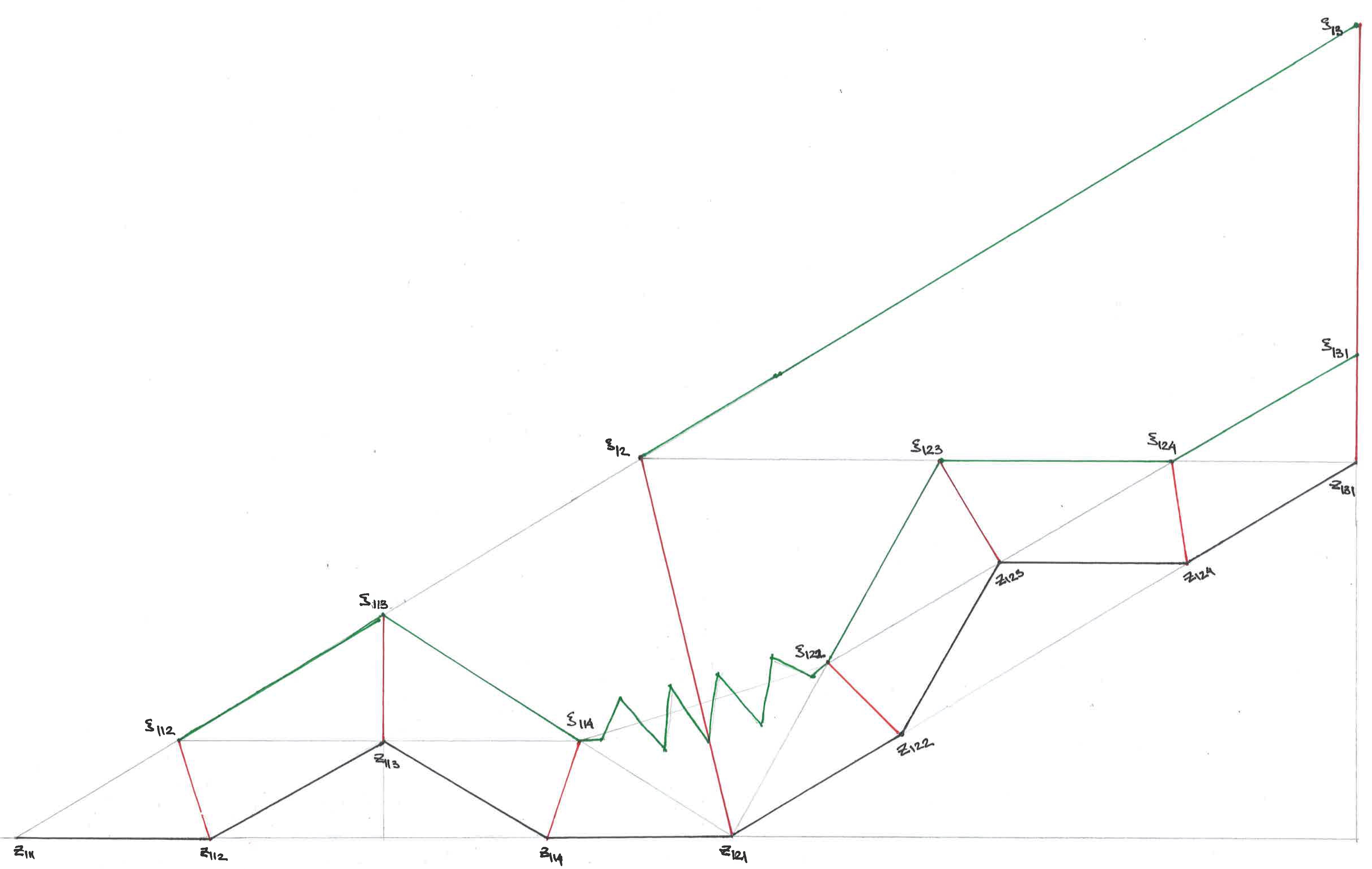}
    \caption{$S^2_\alpha$ and the curve $s_2$}
    \label{fig:snow 1}
\end{figure}

{\bfseries Step 2.} Again, recall that, by $\eqref{snowdef1}$, we have
$$S^2_\alpha = [z_{111}, z_{112}] \cup [z_{112}, z_{113}] \cup [z_{113}, z_{114}] \cup [z_{114}, z_{121}] \cup$$
$$[z_{121}, z_{122}] \cup [z_{122}, z_{123}] \cup [z_{123}, z_{124}] \cup [z_{124}, z_{131}] \cup [z_{131}, z_{132}] \cup [z_{132}, z_{133}] \cup [z_{133}, z_{134}] \cup [z_{134}, z_{141}] \cup$$
$$[z_{141}, z_{142}] \cup [z_{142}, z_{143}] \cup [z_{143}, z_{144}] \cup [z_{144}, z_{211}] \quad \mbox{with}$$
$$z_{111} = z_{11}, \; z_{112} = z_{111} + l(z_{12} - z_{11}), \; z_{113} = z_{112} + l(z_{12} - z_{11})e^{i\alpha}, \; z_{114} = z_{113} + l(z_{12} - z_{11})e^{-i\alpha},$$
$$z_{121} = z_{12}, \; z_{122} = z_{121} + l(z_{13} - z_{12}), \; z_{123} = z_{122} + l(z_{13} - z_{12})e^{i\alpha}, \; z_{124} = z_{123} + l(z_{13} - z_{12})e^{-i\alpha},$$
$$z_{131} = z_{13}, \; z_{141} = z_{14}, \; \mbox{and} \quad z_{13i} = -\overline{z_{12(6 - i)}}, \; z_{14i} = -\overline{z_{11(6 - i)}} \quad \mbox{for} \; i = 2, 3, 4.$$
We now start defining $s_2$. Our (vague) goal is to have $s_k$ sort of more and more self-similar, just as the iterations $S^k_\alpha$. This will be achieved quasi-automatically on step 3 and later, as announced, but for now we need to build $s_2$ by hand and we still want to have some similarity to $s_1$. So, since the piece of $S^2_\alpha$ we got above the segment $[z_{11}, z_{12}]$ is similar to $S^1_\alpha$ above $[z_1, z_2]$, we want to have in $s_2$ two segments $[\zeta_{112}, \zeta_{113}]$ and $[\zeta_{113}, \zeta_{114}]$ with 
$$\zeta_{112} = z_{111} + l(z_{12} - z_{11})e^{i\alpha}, \; \zeta_{113} = \zeta_{112} + \frac{l}{\cos{\alpha}}(z_{12} - z_{11})e^{i\alpha}, \; \zeta_{114} = \zeta_{113} + \frac{l}{\cos{\alpha}}(z_{12} - z_{11})e^{-i\alpha}.$$
The piece of $S^2_\alpha$ above the segment $[z_{12}, z_{13}]$ is also similar to $S^1_\alpha$, so we also want to have in $s_2$ two segments $[\zeta_{122}, \zeta_{123}]$ and $[\zeta_{123}, \zeta_{124}]$ with
$$\zeta_{122} = z_{121} + l(z_{13} - z_{12})e^{i\alpha}, \; \zeta_{123} = \zeta_{122} + \frac{l}{\cos{\alpha}}(z_{13} - z_{12})e^{i\alpha}, \; \zeta_{124} = \zeta_{123} + \frac{l}{\cos{\alpha}}(z_{13} - z_{12})e^{-i\alpha}.$$
Analogously, we add to $s_2$ segments $[\zeta_{132}, \zeta_{133}]$, $[\zeta_{133}, \zeta_{134}]$, $[\zeta_{142}, \zeta_{143}]$ and $[\zeta_{143}, \zeta_{144}]$ with $\zeta_{13i} = -\overline{\zeta_{12(6 - i)}}$ and $\zeta_{14i} = -\overline{\zeta_{11(6 - i)}}$ for $1 = 2, 3, 4$. Go to figure $\ref{fig:snow 1}$ below to see how the defined pieces to the left of $\{x = 0\}$ look like.

We finish the definition of $s_2$ by connecting $\zeta_{114}$ with $\zeta_{122}$, $\zeta_{134}$ with $\zeta_{142}$, and $\zeta_{124}$ and $\zeta_{132}$ with the symmetry axis $\{x = 0\}$. The last one is easy: we add to $s_2$ segments $[\zeta_{124}, \zeta_{131}]$ and $[\zeta_{131}, \zeta_{132}]$ with $\zeta_{131} = \zeta_{124} + \frac{l}{\cos{\alpha}}(z_{13} - z_{12})e^{i\alpha}$. With other parts we need to be more careful. To start with, of course, we want the two curves $(\zeta_{134}, \zeta_{142})$ and $(\zeta_{114}, \zeta_{122})$ to be symmetric to each other (with respect to $\{x = 0\}$): $(\zeta_{134}, \zeta_{142}) = -\overline{(\zeta_{114}, \zeta_{122})}$, so we concentrate on $(\zeta_{114}, \zeta_{122})$. We need the curve $(\zeta_{114}, \zeta_{122})$ to satisfy the following properties:
\begin{enumerate}
    \item[a)] first, symmetry with respect to the bisection of $[\zeta_{114}, \zeta_{122}]$ (which is the line $(z_{12}, \zeta_{12})$),
    \item[b)] second, to have the same length as the piece of $s_2$ composed of $[\zeta_{121}, \zeta_{113}]$ and $[\zeta_{113}, \zeta_{114}]$ together,
    \item[c)] and third, in small balls centred at $\zeta_{114}$ and $\zeta_{122}$, where small means that the radius is much less than $l^2$, the curve is symmetric with respect to the axis $(z_{114}, \zeta_{114})$ or $(z_{122}, \zeta_{122})$ to the segment $[\zeta_{114}, \zeta_{113}]$ or $[\zeta_{122}, \zeta_{123}]$ respectively. 
\end{enumerate}
The purpose of property a) is to help building the intermediate curves between $s_2$ and $s_3$ later on. Intuition behind property b) is that it makes a Brownian traveller modulated by our elliptic operator land with the same probability on $S_\alpha$ everywhere. Property c) we need in order to preserve all the different symmetries of the system after we mollify the construction in Section 5. We will comment more on all of them later (so as not to overcrowd the part where we define everything). We also do not want $s_2$ to intersect with all the other green curves of the $\{s_k\}$ family, which will imply another (not too restrictive) condition we formulate precisely at the next step; vaguely, this condition says that the curve stays close to the segment $[\zeta_{114}, \zeta_{122}]$ and does not wiggle around too much. This implies in particular that this part of $s_2$ does not intersect the lines $z_{121} - te^{-i\alpha}$, $z_{121} + te^{i2\alpha}$ and $z_{111} + te^{i\alpha}$, $t \in \R$. The length of $[\zeta_{121}, \zeta_{113}]$ and $[\zeta_{113}, \zeta_{114}]$ is $\frac{l^2}{\cos{\alpha}}$. According to our restrictions above, we can choose any curve of the length $\frac{l^2}{\cos{\alpha}}$ which begins with the line segment according to c), does not intersect all the lines mentioned (plus behaves reasonably as indicated above) and connects $\zeta_{114}$ and the bisection of $[\zeta_{114}, \zeta_{122}]$, and then reflect it with respect to the latter bisector to get a curve which connects $\zeta_{114}$ with $\zeta_{122}$. Note that $|[z_{121}, \zeta_{114}]| = l^2 < \frac{l^2}{\cos{\alpha}}$, so whatever our curve is, it cannot be a straight line. On the picture below we chose a jagged broken line, just because it is the simplest choice in terms of geometry; but we repeat that we have some degree of freedom here. Again, see figure $\ref{fig:snow 1}$ to have an idea about how a variant of the curve $(\zeta_{114}, \zeta_{122})$ looks like. This completes the definition of $s_2$ in the considered region. 

At the end of the step 2 we define 12 new red lines $\theta_{112} = [z_{112}, \zeta_{112}]$, $\theta_{113} = [z_{113}, \zeta_{113}]$, $\theta_{114} = [z_{114}, \zeta_{114}]$, $\theta_{122} = [z_{122}, \zeta_{122}]$, $\theta_{123} = [z_{123}, \zeta_{123}]$, $\theta_{124} = [z_{124}, \zeta_{124}]$, and finally $\theta_{13i} = -\overline{\theta_{12(6 - i)}}$, $\theta_{14i} = -\overline{\theta_{11(6 - i)}}$ for $i = 2, 3, 4$. Figure $\ref{fig:snow 1}$ is the sketch of the part of what we have constructed thus far which lies to the left of $\{x = 0\}$. To have the full picture it suffices to reflect the one below with respect to $\{x = 0\}$.

Before we move on to Step 3, observe that
\begin{equation}\label{snowdef2}
    S^2_\alpha = F(S^1_\alpha, z_{11}),
\end{equation}
and we can even define all the points of the second generation on $S^2_\alpha$ using $\eqref{snowdef2}$ as a definition: 
$z_{1ij} = z_{F_i(S^1_\alpha)}$ for $z = z_{1j}$, $i, j = 1, 2, 3, 4$ (see Definition $\ref{F}$). 

\begin{figure}[h]
    \centering
    \includegraphics[scale = 0.08]{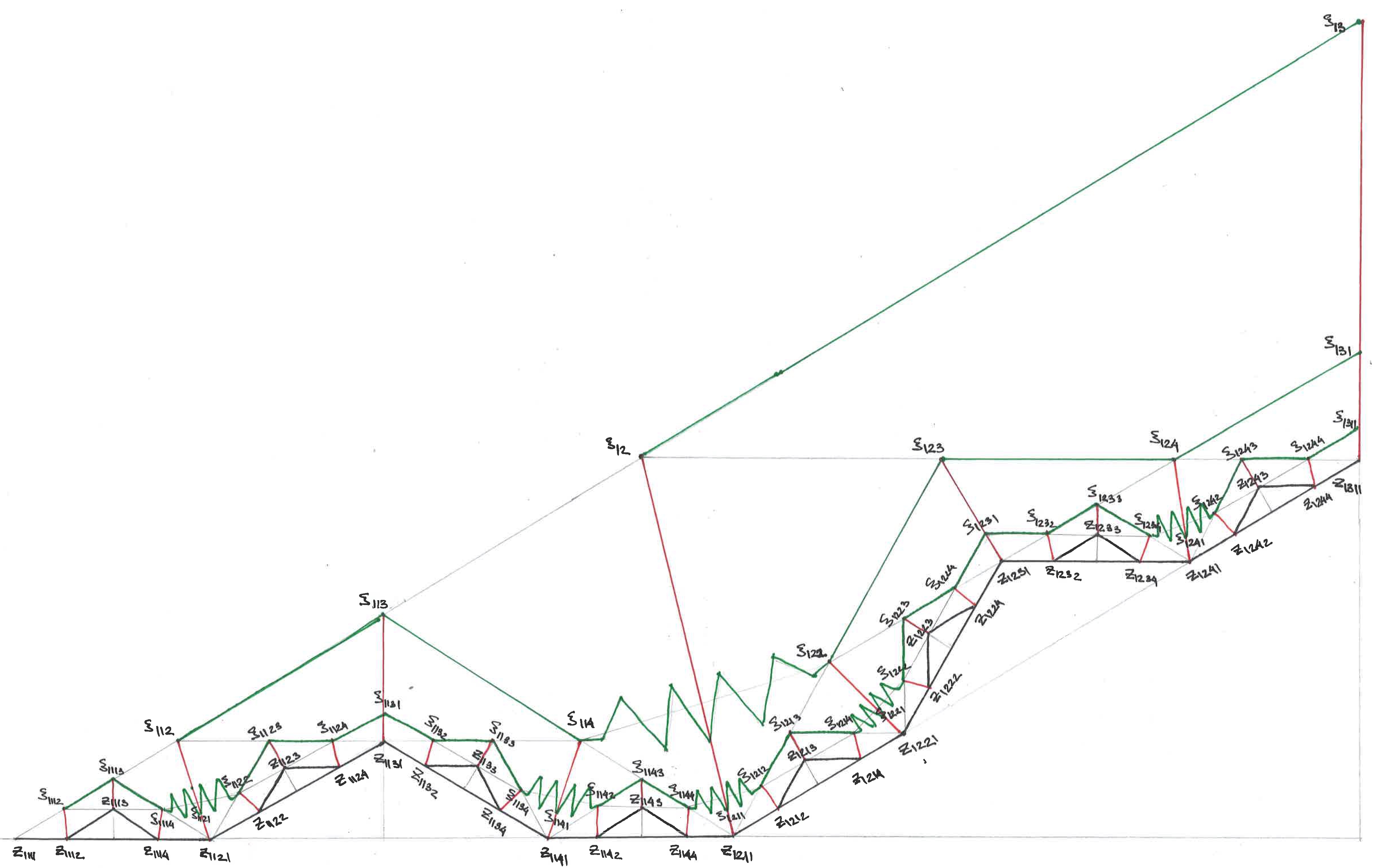}
    \caption{$S^3_\alpha$ and the curve $s_3$}
    \label{fig:snow 2}
\end{figure}

{\bfseries Step 3.} From this step on, we will exploit the transformation $F$ a lot. We define (a piece of) the new green curve
$$s_3 = F(s_2, \zeta_{1112}) \quad \mbox{with} \quad \zeta_{1112} = z_{1111} + l(z_{112} - z_{111})e^{i\alpha}, \quad \mbox{and}$$
\begin{equation}\label{s3}
    \zeta_{w_1iw_2w_3} = z_{F_i(s_2)} \quad \mbox{for} \; z = \zeta_{w_1w_2w_3}, i = 1, 2, 3, 4; \; w = w_1w_2w_3 \in \mathcal{I}_2, \; w_1 = 1.
\end{equation}
We also define new red curves
\begin{equation}
    \theta_{w_1iw_2w_3} = [z_{w_1iw_2w_3}, \zeta_{w_1iw_2w_3}] \quad \mbox{for} \; z = \zeta_{w}, i = 1, 2, 3, 4; \; w \in \mathcal{I}_2, \; w_1 = 1.
\end{equation}
The curve $s_3$ is however not complete, because by the definition $\eqref{s3}$ it is not yet connected: we need to join $\zeta_{i - 1}$ and $\zeta_{i + 1}$ for $i = (1211), (1411)$ and $(1311)$. We do it in the following way:
$$(\zeta_{1244}, \zeta_{1312}) = ((\zeta_{124}, \zeta_{132}) - z_{131})l + z_{1311},$$
\begin{equation}
    (\zeta_{1144}, \zeta_{1212}) = ((\zeta_{114}, \zeta_{122}) - z_{121})l + z_{1211}, \; \mbox{and} \; (\zeta_{1344}, \zeta_{1412}) = -\overline{(\zeta_{1144}, \zeta_{1212})},
\end{equation}
where $(\zeta_{114}, \zeta_{122})$ and $(\zeta_{124}, \zeta_{132})$ are pieces of $s_2$ between the two points indicated. See figure $\ref{fig:snow 2}$ for a sketch of $S^3_\alpha$, the curves $s_1, s_2$ and $s_3$ and all the red curves defined up to that point to the left of the axis $\{x = 0\}$.

We now can complete the set of exact conditions we impose on the curve $(\zeta_{114}, \zeta_{122})$. At Step 2 we described vaguely the condition we are about to specify as ``this part of the curve $s_2$ does not intersect the other curves $s_k$''. In particular, $(\zeta_{114}, \zeta_{122})$ should not intersect our newly defined $(\zeta_{1144}, \zeta_{1212})$. Suppose that we can treat the curve $(\zeta_{114}, \zeta_{122})$ as a graph of a function $f_1$ over the segment $[\zeta_{114}, \zeta_{122}]$; therefore, by similarity, we can treat $(\zeta_{1144}, \zeta_{1212})$ as a graph of a function $f_2$ over the segment $[\zeta_{1144}, \zeta_{1212}]$. Then $\max{f_2} = l\max{f_1}$ and $\min{f_2} = l\min{f_1}$. We ask in addition $\max{f_1} = |\min{f_1}|$ so not to complicate things. It is clear from the picture that, for the two graphs not to intersect, we need to impose the condition 
$$|\min{f_1}| + \max{f_2} = \max{f_1}(1 + l) < \dist([\zeta_{114}, \zeta_{122}], [\zeta_{1144}, \zeta_{1212}]) = $$ $$ = \cos{\left(\frac{\pi}{2} - \frac{3\alpha}{2}\right)}\left(|z_{1211} - \zeta_{114}| - |z_{1211} - \zeta_{1144}|\right) = \sin{(3\alpha/2)}(l^2\cos{\alpha} - l^3\cos{\alpha}) = $$ $$ = l^2(1 - l)\cos{\alpha}\sin{(3\alpha/2)},$$
which implies
\begin{equation}\label{cond_graph}
    \max{f_1} = |\min{f_1}| < l^2 \frac{1 - l}{1 + l}\cos{\alpha}\sin{(3\alpha/2)}.
\end{equation}
To summarise, the exact condition we impose on the curve $(\zeta_{114}, \zeta_{122})$, in addition to the symmetry and the length conditions, is that $(\zeta_{114}, \zeta_{122})$ is a graph of a (continuous and smooth except for a finite number of points) function $f_1$ and that $f_1$ admits $\eqref{cond_graph}$. To complete the list above, we will call it condition d).

\medskip

Again, observe that 
\begin{equation}
    S^3_\alpha = F(S^2_\alpha, z_{111}),
\end{equation}
and we could have even defined all the points of the third generation on $S^3_\alpha$ using this relation as a definition:
\begin{equation}
    z_{w_1iw_2w_3} = z_{F_i(S^2_\alpha, z_{111})} \quad \mbox{for} \; z = z_{w_1w_2w_3}, i = 1, 2, 3, 4; \; w = w_1 w_2 w_3 \in \mathcal{I}_2, \; w_1 = 1.
\end{equation}

We pass on to the general case. {\bfseries Step k.} Define 
$$s_k = F(s_{k - 1}, \zeta_{1 \dots 12}) \quad \mbox{with} \quad \zeta_{1 \dots 12} = z_{1 \dots 1} + l(z_{1 \dots 2} - z_{1\dots 1})e^{i\alpha},$$
where $\dots$ denotes $k - 2$ ones in a row, and
\begin{equation}\label{s_k}
    \zeta_{w_1iw_2\dots w_k} = z_{F_i(s_{k - 1})} \quad \mbox{for} \; z = \zeta_{w_1 \dots w_k}, i = 1, 2, 3, 4; \; w = w_1\dots w_k \in \mathcal{I}_{k - 1}, \; w_1 = 1.
\end{equation}
Define new red curves
\begin{equation}\label{theta_k}
    \theta_{w_1iw_2 \dots w_k} = [z_{w_1iw_2 \dots w_k}, \zeta_{w_1iw_2 \dots w_k}] \quad \mbox{for} \; z = \zeta_{w}, i = 1, 2, 3, 4; \; w \in \mathcal{I}_{k - 1}, \; w_1 = 1.
\end{equation}
To complete the definition of $s_k$ we need to draw curves between $\zeta_{i - 1}$ and $\zeta_{i + 1}$ for \; $i = 121\dots, \; 141\dots$ and $131\dots$, where $\dots$ once again denotes $k - 2$ ones in a row. We draw
$$(\zeta_{i - 1}, \zeta_{i + 1}) = ((\zeta_{124 \dots 4}, \zeta_{13 \dots 12}) - z_{131 \dots})l + z_{131 \dots} \quad \mbox{for} \;  i = 131\dots,$$
$$(\zeta_{i - 1}, \zeta_{i + 1}) = ((\zeta_{114 \dots 4}, \zeta_{12 \dots 12}) - z_{121 \dots})l + z_{121 \dots} \quad \mbox{for} \;  i = 121\dots,$$
where the first $\dots$ denote $k - 3$ fours in a row, and all the others -- $k - 3$ ones, and
\begin{equation}\label{complete_s_k}
     (\zeta_{i - 1}, \zeta_{i + 1}) = -\overline{(\zeta_{1144 \dots 4}, \zeta_{121 \dots 12})} \quad \mbox{for} \; i = 141\dots.
\end{equation}
We finish with the observation that
\begin{equation}\label{S^k_alpha}
    S^k_\alpha = F(S^{k - 1}_\alpha, z_1),
\end{equation}
and 
\begin{equation}\label{sk}
    z_{w_1iw_2 \dots w_k} = z_{F_i(S^{k - 1}_\alpha, z_1)} \quad \mbox{for} \; z = z_{w_1w_2 \dots w_k}, i = 1, 2, 3, 4; \; w = w_1w_2\dots w_k \in \mathcal{I}_{k - 1}, \; w_1 = 1.
\end{equation}

\section{Tiling of the space}

In this section we prepare the terrain for completing the net of curves in a neighbourhood of the part of the fractal generated by the unit interval $[z_1, z_2]$. Do to this, we divide certain regions composed from the cells of our large-scale net from the previous section (formed by adjacent green curves $s_k$ and $s_{k + 1}$ and red curves ${\theta_w}$) into two types. All the regions of the same type will be similar to each other. By similar we mean the usual thing in Euclidean geometry: $A$ and $B$ are similar if one is an image of the other under a composition $h \circ m$, where $h$ is a scaling transformation and $m$ -- a motion. Once this is done, we will explain in Section 6 how to fill in a region of each of those two types with green and red curves, so that in addition those intermediate curves glue nicely to the surrounding picture. It will suffice to give this explanation only once for each type of region, since we can complete the net in the rest of the regions of the same type by the similarity and using our transformation $F$. Thus we will cover with the net all the target neighbourhood of $[z_1, z_2]$. We will call ``tile'' a representative of one type of regions, because the two types of regions tile the neighbourhood of $\partial \Omega$ we are interested in.

\begin{figure}[h]
    \centering
    \includegraphics[scale = 0.2]{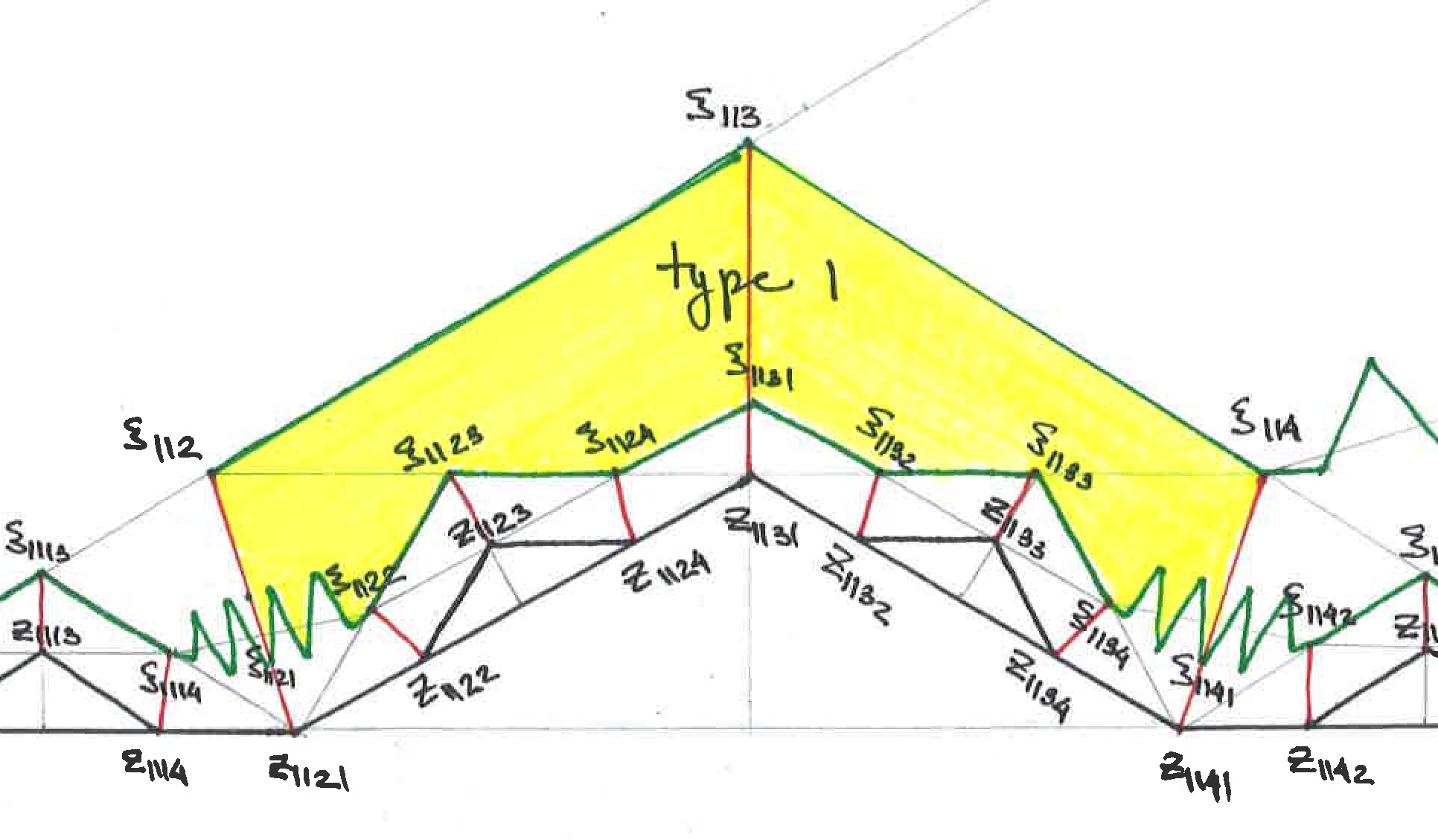}
    \caption{Tile T1}
    \label{fig:tile type 1}
\end{figure}

{\bfseries Type 1 tiles.} We will call tiles of type 1 all the regions which are similar to the figure 
\begin{equation}
    T1 = \zeta_{1121}\zeta_{112}\zeta_{113}\zeta_{114}\zeta_{1141}\zeta_{1134}\zeta_{1133}\zeta_{1132}\zeta_{1131}\zeta_{1124}\zeta_{1123}\zeta_{1122},
\end{equation}
composed by two red segments $[\zeta_{1121}\zeta_{112}]$ and $[\zeta_{114}, \zeta_{1141}]$, and two pieces of green curves $s_2$ and $s_3$: $[\zeta_{112}, \zeta_{113}] \cup [\zeta_{113}, \zeta_{114}]$ and $(\zeta_{1121}, \zeta_{1122}, \zeta_{1123}, \zeta_{1124}, \zeta_{1131}) \cup (\zeta_{1131}, \zeta_{1132}, \zeta_{1133}, \zeta_{1134}, \zeta_{1141})$. Note that $T1$ is symmetric with respect to the axis $(\zeta_{113}, \zeta_{1131})$. See figure $\ref{fig:tile type 1}$.

In the neighbourhood of the snowflake which we consider, between the curves $s_2$ and $s_3$, there are other four tiles of type 1 of the same size (similar to $T1$ with coefficient 1):
$$\zeta_{1221}\zeta_{122}\zeta_{123}\zeta_{124}\zeta_{1241}\zeta_{1234}\zeta_{1233}\zeta_{1232}\zeta_{1231}\zeta_{1224}\zeta_{1223}\zeta_{1222} = (T1 - \zeta_{112})e^{i\alpha} + \zeta_{122},$$
$$\zeta_{1321}\zeta_{132}\zeta_{133}\zeta_{134}\zeta_{1341}\zeta_{1334}\zeta_{1333}\zeta_{1332}\zeta_{1331}\zeta_{1324}\zeta_{1323}\zeta_{1322} = (T1 - \zeta_{112})e^{-i\alpha} + \zeta_{132},$$
$$\zeta_{1421}\zeta_{142}\zeta_{143}\zeta_{144}\zeta_{1441}\zeta_{1434}\zeta_{1433}\zeta_{1432}\zeta_{1431}\zeta_{1424}\zeta_{1423}\zeta_{1422} = (T1 - \zeta_{112}) + \zeta_{142}, \quad \mbox{and}$$
$$\zeta_{1241}\zeta_{124}\zeta_{131}\zeta_{132}\zeta_{1321}\zeta_{1314}\zeta_{1313}\zeta_{1312}\zeta_{1311}\zeta_{1244}\zeta_{1243}\zeta_{1242} = (T1 - \zeta_{112}) + \zeta_{124}.$$

\medskip

\begin{figure}[h]
    \centering
    \includegraphics[scale = 0.38]{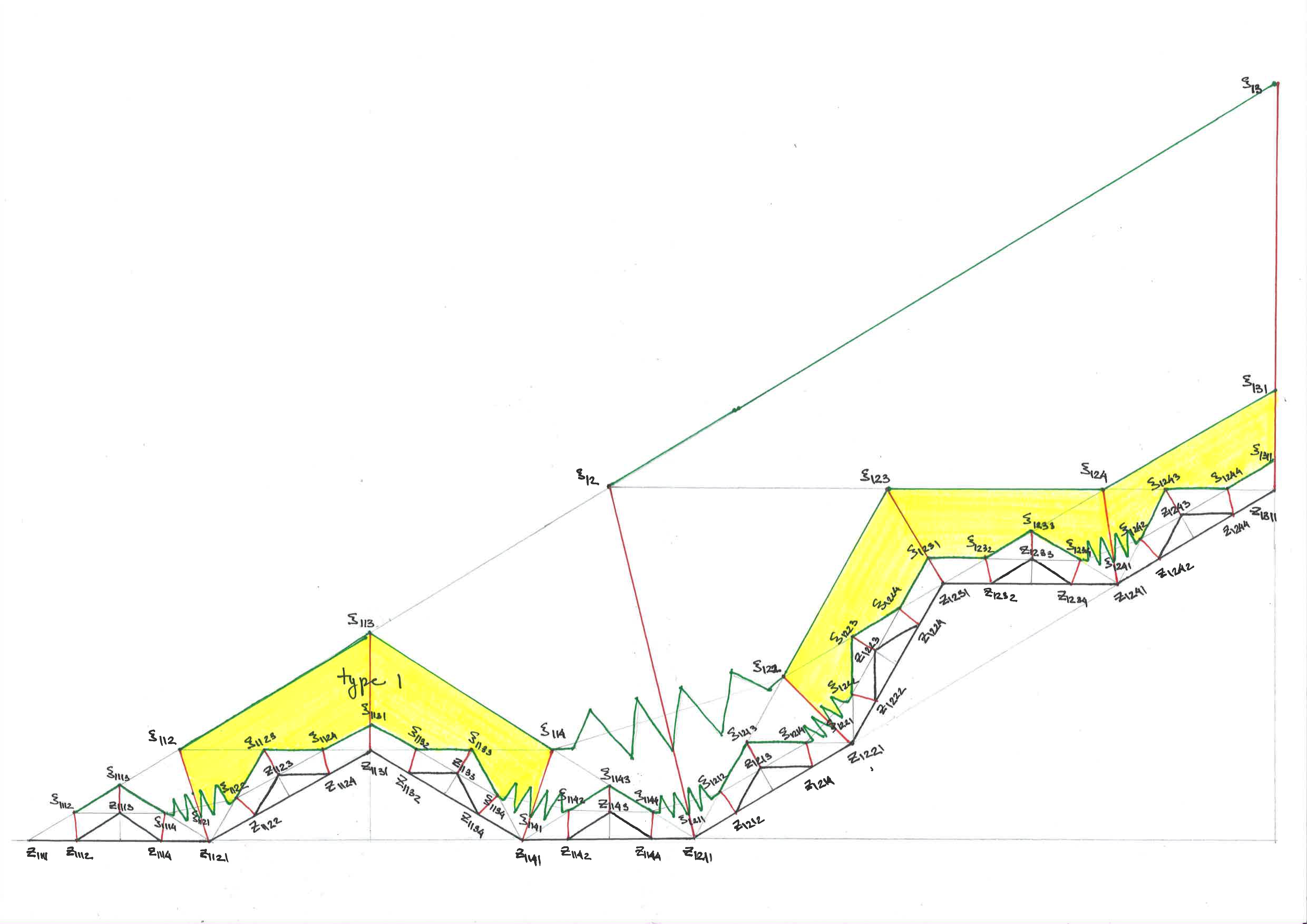}
    \caption{Tiles of type 1 between the curves $s_2$ and $s_3$}
    \label{fig:type 1 s2_s3}
\end{figure}

Figure $\ref{fig:type 1 s2_s3}$ shows all the tiles of type 1 in the stripe between $s_2$ and $s_3$ to the left of the axis $\{x = 0\}$.

\begin{figure}[h]
    \centering
    \includegraphics[scale = 0.2]{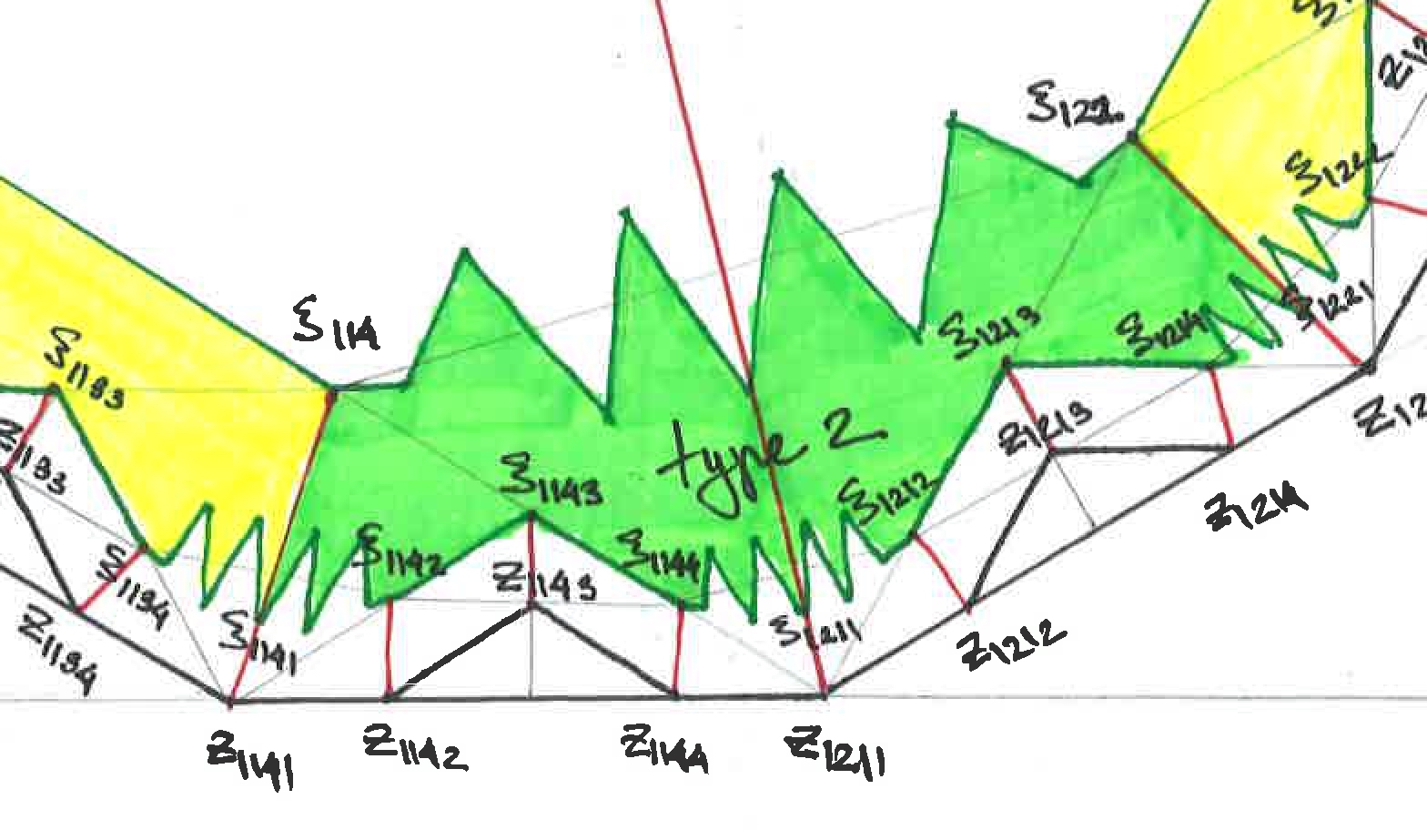}
    \caption{Tile T2}
    \label{fig:tile type 2}
\end{figure}

{\bfseries Type 2 tiles.} We will call tiles of type 2 all the regions which are similar to the figure 
\begin{equation}
    T2 = \zeta_{1141}\zeta_{114}\zeta_{121}\zeta_{122}\zeta_{1221}\zeta_{1214}\zeta_{1213}\zeta_{1212}\zeta_{1211}\zeta_{1144}\zeta_{1143}\zeta_{1142},
\end{equation}
composed by two red segments $[\zeta_{1141}\zeta_{114}]$ and $[\zeta_{122}\zeta_{1221}]$, and two pieces of green curves $s_2$ and $s_3$: the curve $[\zeta_{114}, \zeta_{121}] \cup [\zeta_{121}, \zeta_{122}]$ and the curve $(\zeta_{1141}, \zeta_{1142}, \zeta_{1143}, \zeta_{1144}, \zeta_{1211}) \cup (\zeta_{1211}, \zeta_{1212}, \zeta_{1213}, \zeta_{1214}, \zeta_{1221})$. Note that $T2$ is symmetric with respect to the axis $(\zeta_{121}, \zeta_{1211})$. See figure $\ref{fig:tile type 2}$.

In the same neighbourhood, between the curves $s_2$ and $s_3$, we have another tile of type 2: 
$$\zeta_{1341}\zeta_{134}\zeta_{141}\zeta_{142}\zeta_{1421}\zeta_{1414}\zeta_{1413}\zeta_{1412}\zeta_{1411}\zeta_{1344}\zeta_{1343}\zeta_{1342} = -\overline{(T2 - z_{1211})} + z_{1411}.$$

\begin{figure}[h]
    \centering
    \includegraphics[scale = 0.38]{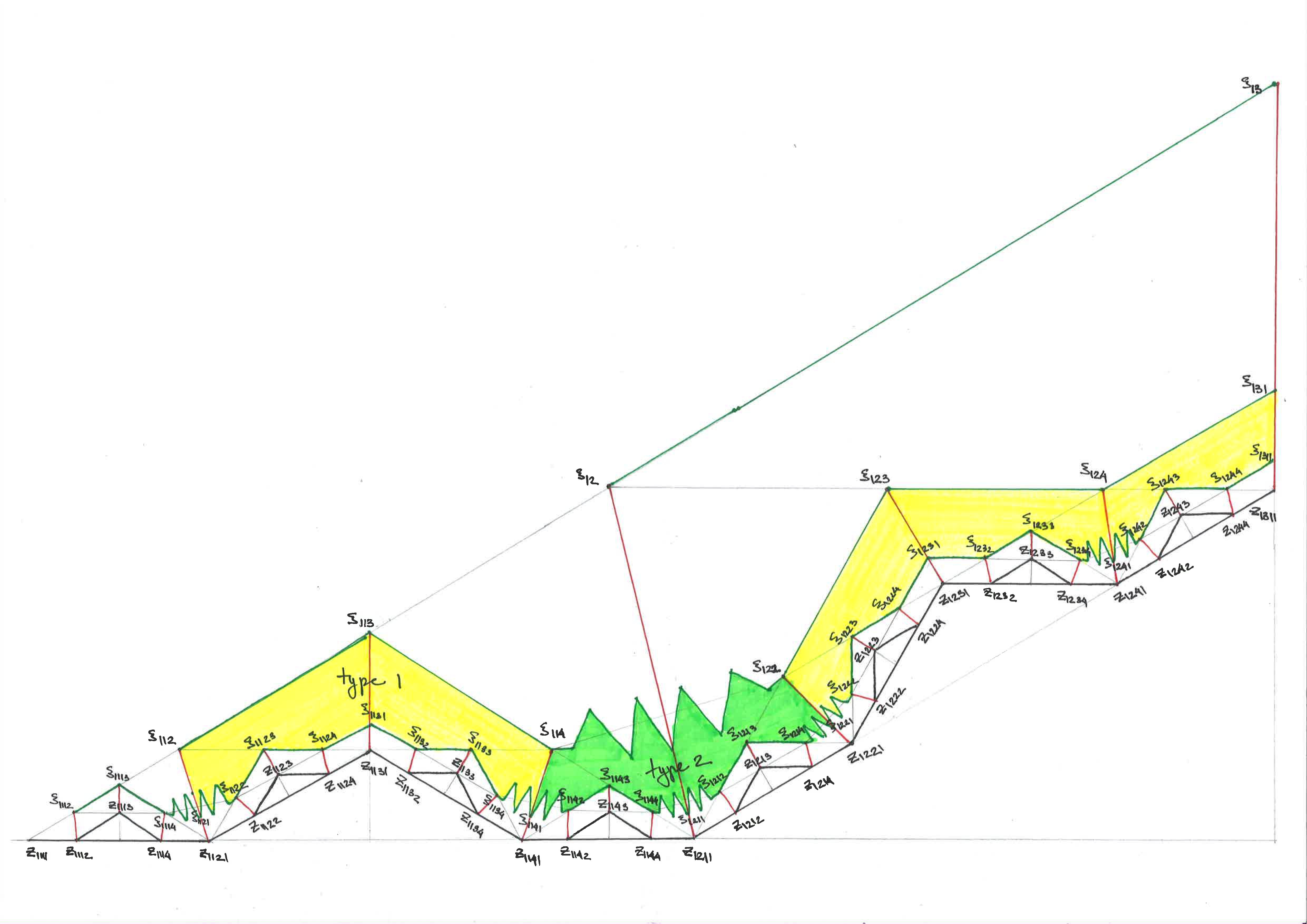}
    \caption{Tiling of the region between the curves $s_2$ and $s_3$}
    \label{fig:tiling s2_s3}
\end{figure}

Note that the five tiles of type 1 and two tiles of type 2 enlisted above tile all the space between $s_2$ and $s_3$ in the neighbourhood of the snowflake restricted by the left half of the triangle $z_1z_2\zeta_{13}$, except for two small regions near $[z_1, \zeta_{13}]$ and $[\zeta_{13}, z_2]$: $\zeta_{1114}\zeta_{1113}\zeta_{112}\zeta_{1121}$, bounded by two red segments, a piece of $s_3$ and a segment on the line $(z_1 \zeta_{13})$, and the region $\zeta_{1433}\zeta_{1434}\zeta_{144}\zeta_{1441} = -\overline{\zeta_{1114}\zeta_{1113}\zeta_{112}\zeta_{1121}}$. So we may assume that for now the tiling of the space between $s_2$ and $s_3$ is complete (modulo this little area just mentioned, with which we deal with later). Figure $\ref{fig:tiling s2_s3}$ shows how the tiling looks like to the left of the axis $\{x = 0\}$.

We now explain how to tile the space between the curves $s_{k - 1}$ and $s_k$, $k > 3$. Once again, we will tile the whole stripe inside the triangle $z_1z_2\zeta_{13}$, except for the regions $\zeta_{1\dots 14}\zeta_{1\dots 13}\zeta_{1 \dots 12}\zeta_{1\dots 21}$ (lying between the line $(z_1, \zeta_{13}$ and the red curve $[\zeta_{1\dots 12}, z_{1\dots 12}]$) and the one centrally symmetric to it with respect to $\{x = 0\}$, $-\overline{\zeta_{1\dots 14}\zeta_{1\dots 13}\zeta_{1 \dots 12}\zeta_{1\dots 21}}$ (lying between the line $(\zeta_{13}, z_2)$ and the red curve $[\zeta_{14\dots 4}, z_{14\dots 4}]$). We do that again by induction, exploiting the definition $\eqref{s_k}$ of $s_k$, that is, the transformation $F$. The tiling of the space between $s_2$ and $s_3$ we have just done is the base.

From now on we will call tiles $T1$ or $T2$ all the figures similar to $T1$ or $T2$ respectively. In what follows we also denote by $St(s_k, s_{k + 1}, I, J)$ the piece of stripe lying between the curves $s_k$ and $s_{k + 1}$ bounded by the red segments $I$ and $J$. 

\begin{figure}[h]
    \centering
    \includegraphics[scale = 0.08]{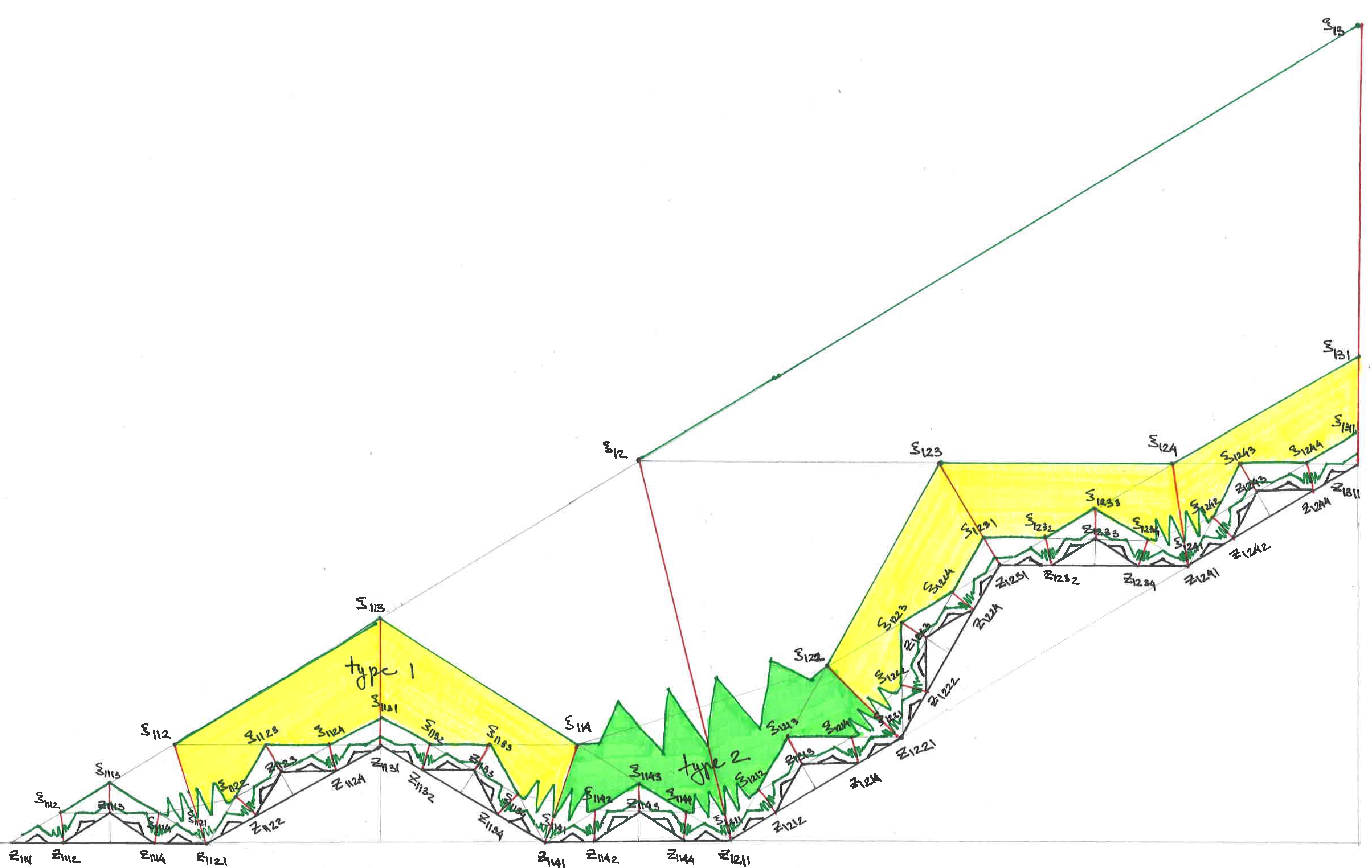}
    \caption{Tiling between the curves $s_2$ and $s_3$, the curve $s_4$}
    \label{fig:tiling s2_s3 plus s4}
\end{figure}

We illustrate what happens on the stripe between $s_3$ and $s_4$ (which corresponds to the first step of induction). Figure $\ref{fig:tiling s2_s3 plus s4}$ represents how the whole picture of the tiling we did on the base step along with the curve $s_4$ looks like.

Observe that, according to $\eqref{s_k}$ and $\eqref{theta_k}$, we have 
$$St(s_3, s_4, [z_{1112}, \zeta_{1112}], [z_{1144}, \zeta_{1144}]) = F_1(St(s_2, s_3, [z_{112}, \zeta_{112}], [z_{144}, \zeta_{144}])),$$
$$St(s_3, s_4, [z_{1212}, \zeta_{1212}], [z_{1244}, \zeta_{1244}]) = F_2(St(s_2, s_3, [z_{112}, \zeta_{112}], [z_{144}, \zeta_{144}])),$$
$$St(s_3, s_4, [z_{1312}, \zeta_{1312}], [z_{1344}, \zeta_{1344}]) = F_3(St(s_2, s_3, [z_{112}, \zeta_{112}], [z_{144}, \zeta_{144}])), \quad \mbox{and}$$
\begin{equation}
    St(s_3, s_4, [z_{1412}, \zeta_{1412}], [z_{1444}, \zeta_{1444}]) = F_4(St(s_2, s_3, [z_{112}, \zeta_{112}], [z_{144}, \zeta_{144}]))
\end{equation}
with $p = z_1$ in the definition of transformation $F$. Therefore, since the figure 
$$St(s_2, s_3, [z_{112}, \zeta_{112}], [z_{144}, \zeta_{144}])$$
is tiled completely by $T_1$ and $T_2$, the whole stripe between $s_3$ and $s_4$ is also almost tiled, except for three figures
$$F_1 = \zeta_{11441}\zeta_{1144}\zeta_{1211}\zeta_{1212}\zeta_{12121}\zeta_{12114}\zeta_{12113}\zeta_{12112}\zeta_{12111}\zeta_{11444}\zeta_{11443}\zeta_{11442},$$
$$F_2 = \zeta_{13441}\zeta_{1344}\zeta_{1411}\zeta_{1412}\zeta_{14121}\zeta_{14114}\zeta_{14113}\zeta_{14112}\zeta_{14111}\zeta_{13444}\zeta_{13443}\zeta_{13442} \quad \mbox{and}$$
$$F_3 = \zeta_{12441}\zeta_{1244}\zeta_{1311}\zeta_{1312}\zeta_{13121}\zeta_{13114}\zeta_{13113}\zeta_{13112}\zeta_{13111}\zeta_{12444}\zeta_{12443}\zeta_{12442}.$$

For those figures, the following identities hold:
\begin{equation}\label{F3}
    F_3 = l(T1 - z_{113}) + z_{1311},
\end{equation}
\begin{equation}\label{F12}
    F_1 = l(T2 - z_{121}) + z_{1211} \quad \mbox{and} \quad F_2 = -\overline{F_1}.
\end{equation}
We explain first $\eqref{F12}$. The second equality there is just a corollary of the symmetry of the entire construction with respect to the axis $\{x = 0\}$. The first one is the combination of the following observations. First, by $\eqref{complete_s_k}$, we have
$$(\zeta_{1144}, \zeta_{1212}) = ((\zeta_{114}, \zeta_{122}) - z_{121})l + z_{1211}, \quad (\zeta_{11444}, \zeta_{12112}) = ((\zeta_{1144}, \zeta_{1212}) - z_{1211})l + z_{12111}.$$
Second, by $\eqref{s_k}$, 
$$(\zeta_{11441}, \zeta_{11442}) \cup [\zeta_{11442}, \zeta_{11443}] \cup [\zeta_{11443}, \zeta_{11444}] = l((\zeta_{1441}, \zeta_{1442}) \cup [\zeta_{1442}, \zeta_{1443}] \cup [\zeta_{1443}, \zeta_{1444}] - z_1) + z_1 = $$ $$= l(l((\zeta_{141}, \zeta_{142}) \cup [\zeta_{142}, \zeta_{143}] \cup [\zeta_{143}, \zeta_{144}] - z_1) + z_{14} - z_1) + z_1, \quad \mbox{and}$$
$$(\zeta_{1141}, \zeta_{1142}) \cup [\zeta_{1142}, \zeta_{1143}] \cup [\zeta_{1143}, \zeta_{1144}] = l((\zeta_{141}, \zeta_{142})\cup [\zeta_{142}, \zeta_{143}] \cup [\zeta_{143}, \zeta_{144}] - z_1) + z_{11}, \quad \mbox{so}$$
$$(\zeta_{11441}, \zeta_{11442}) \cup [\zeta_{11442}, \zeta_{11443}] \cup [\zeta_{11443}, \zeta_{11444}] = $$ 
$$l(((\zeta_{1141}, \zeta_{1142}) \cup [\zeta_{1142}, \zeta_{1143}] \cup [\zeta_{1143}, \zeta_{1144}] - z_{11}) + z_{14} - z_1) + z_1 = $$ $$ = l((\zeta_{1141}, \zeta_{1142}) \cup [\zeta_{1142}, \zeta_{1143}] \cup [\zeta_{1143}, \zeta_{1144}]) + l(z_{14} - 2z_1) + z_1.$$
One can compute that $l(z_{14} - 2z_1) + z_1 = -l z_{12} + z_{12}$, so we indeed have
$$(\zeta_{11441}, \zeta_{11442}) \cup  [\zeta_{11442}, \zeta_{11443}] \cup [\zeta_{11443}, \zeta_{11444}] = l((\zeta_{1141}, \zeta_{1142}) \cup [\zeta_{1142}, \zeta_{1143}] \cup [\zeta_{1143}, \zeta_{1144}] - z_{12}) + z_{12}.$$
Third, again by $\eqref{s_k}$, 
$$[\zeta_{12112}, \zeta_{12113}] \cup [\zeta_{12113}, \zeta_{12114}] \cup (\zeta_{12114}, \zeta_{12121}) = l([\zeta_{1112}, \zeta_{1113}] \cup [\zeta_{1113}, \zeta_{1114}] \cup (\zeta_{1114}, \zeta_{1121}) - z_1)e^{i\alpha} + z_{12} = $$ $$= l(l([\zeta_{112}, \zeta_{113}] \cup [\zeta_{113}, \zeta_{114}] \cup (\zeta_{114}, \zeta_{121}) - z_1) + z_1 - z_1)e^{i\alpha} + z_{12} =$$ 
$$ l^2([\zeta_{112}, \zeta_{113}] \cup [\zeta_{113}, \zeta_{114}] \cup (\zeta_{114}, \zeta_{121}) - z_1)e^{i\alpha} + z_{12}, \quad \mbox{and}$$
$$[\zeta_{1212}, \zeta_{1213}] \cup [\zeta_{1213}, \zeta_{1214}] \cup (\zeta_{1214}, \zeta_{1221}) = l([\zeta_{112}, \zeta_{113}] \cup [\zeta_{113}, \zeta_{114}] \cup (\zeta_{114}, \zeta_{121}) - z_1)e^{i\alpha} + z_{12}, \quad \mbox{so}$$
$$[\zeta_{12112}, \zeta_{12113}] \cup [\zeta_{12113}, \zeta_{12114}] \cup (\zeta_{12114}, \zeta_{12121}) = l([\zeta_{1212}, \zeta_{1213}] \cup [\zeta_{1213}, \zeta_{1214}] \cup (\zeta_{1214}, \zeta_{1221}) - z_{12}) + z_{12}.$$
The definition $\eqref{theta_k}$ and analogous computations for the red sides of $F_1$ finally give us the first identity in $\eqref{F12}$. So $F_1$ and $F_2$ are also tiles $T2$. 

\begin{figure}[h]
    \centering
    \includegraphics[scale = 0.08]{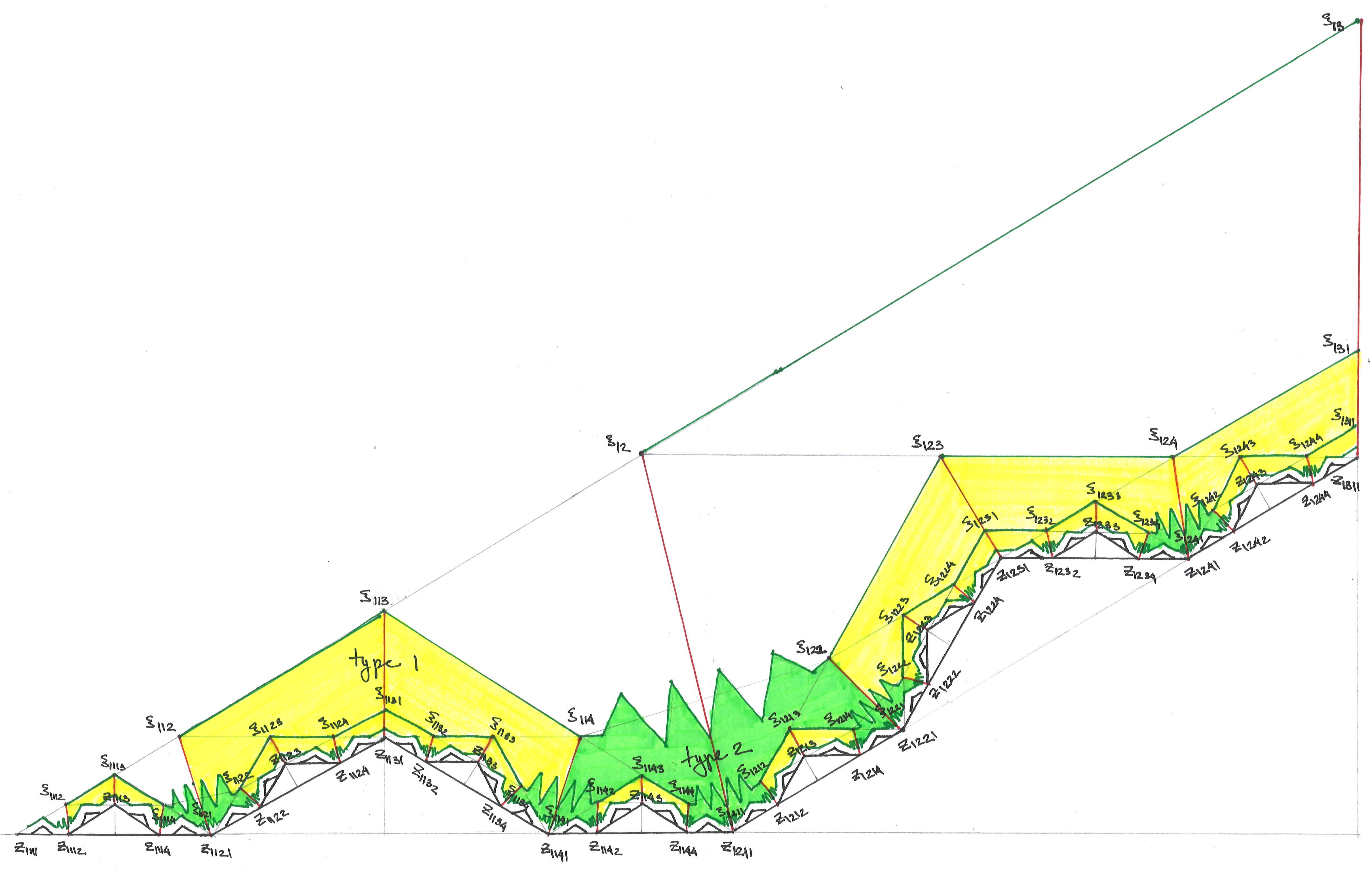}
    \caption{Tiling between the curves $s_2$ and $s_4$}
    \label{fig:tiling s2_s4}
\end{figure}

We now treat $\eqref{F3}$. By $\eqref{complete_s_k}$,
$$[\zeta_{1244}, \zeta_{1311}] \cup [\zeta_{1311}, \zeta_{1312}] = l([\zeta_{124}, \zeta_{131}] \cup [\zeta_{131}, \zeta_{132}] - z_{131}) + z_{1311}, \quad \mbox{and}$$ 
$$ [\zeta_{12444}, \zeta_{13111}] \cup [\zeta_{13111}, \zeta_{13112}] = l([\zeta_{1244}, \zeta_{1311}] \cup [\zeta_{1311}, \zeta_{1312}] - z_{1311}) + z_{13111}.$$
Next, 
$$(\zeta_{12441}, \zeta_{12442}) \cup [\zeta_{12442}, \zeta_{12443}] \cup [\zeta_{12443}, \zeta_{12444}] = l((\zeta_{1441}, \zeta_{1442}) \cup [\zeta_{1442}, \zeta_{1443}] \cup [\zeta_{1443}, \zeta_{1444}] - z_1)e^{i\alpha} + z_{12} $$
$$ = l((l((\zeta_{141}, \zeta_{142}) \cup [\zeta_{142}, \zeta_{143}] \cup [\zeta_{143}, \zeta_{144}] - z_1) + z_{14}) - z_1)e^{i\alpha} + z_{12}.$$
$$l((\zeta_{141}, \zeta_{142}) \cup [\zeta_{142}, \zeta_{143}] \cup [\zeta_{143}, \zeta_{144}] - z_1)e^{i\alpha} + z_{12} = (\zeta_{1241}, \zeta_{1242}) \cup [\zeta_{1242}, \zeta_{1243}] \cup [\zeta_{1243}, \zeta_{1244}].$$
This gives
$$(\zeta_{12441}, \zeta_{12442}) \cup [\zeta_{12442}, \zeta_{12443}] \cup [\zeta_{12443}, \zeta_{12444}] = l((\zeta_{1241}, \zeta_{1242}) \cup [\zeta_{1242}, \zeta_{1243}] \cup [\zeta_{1243}, \zeta_{1244}] - z_{12}) $$ $$ + l(z_{14} - z_1)e^{i\alpha} + z_{12} = l((\zeta_{1241}, \zeta_{1242}) \cup [\zeta_{1242}, \zeta_{1243}] \cup [\zeta_{1243}, \zeta_{1244}] - z_{13}) + $$
$$lz_{13} - lz_{12} + l(z_{14} - z_1)e^{i\alpha} + z_{12}$$
A computation shows that $lz_{13} - lz_{12} + l(z_{14} - z_1)e^{i\alpha} + z_{12} = z_{13}$, so we have
$$(\zeta_{12441}, \zeta_{12442}) \cup [\zeta_{12442}, \zeta_{12443}] \cup [\zeta_{12443}, \zeta_{12444}] = l((\zeta_{1241}, \zeta_{1242}) \cup [\zeta_{1242}, \zeta_{1243}] \cup [\zeta_{1243}, \zeta_{1244}] - z_{13}) + z_{13}.$$
Similar computations show that
$$[\zeta_{13112}, \zeta_{13113]}] \cup [\zeta_{13113}, \zeta_{13114}] \cup (\zeta_{13114}, \zeta_{13121}) = l([\zeta_{1312}, \zeta_{1313}] \cup [\zeta_{1313}, \zeta_{1314}] \cup (\zeta_{1314}, \zeta_{1321}) - z_{13}) + z_{13},$$
and the same holds for the red sides of $F_2$, which combined gives $\eqref{F3}$ and that $F_3$ is a $T1$ tile. This completes the tiling of the stripe between $s_3$ and $s_4$. Figure $\ref{fig:tiling s2_s4}$ shows how the picture of this tiling looks like to the left of $\{x = 0\}$.

The proper induction step does not differ from what we just did for $St(s_3, s_4, [z_{1112}, \zeta_{1112}], [\zeta_{1444}, \zeta_{1444}])$. According to $\eqref{s_k}$ and $\eqref{theta_k}$, we have 
$$St(s_k, s_{k + 1}, [z_{11\dots2}, \zeta_{11\dots2}], [z_{11\dots 4}, \zeta_{11\dots 4}]) = F_1(St(s_{k - 1}, s_k, [z_{1\dots 2}, \zeta_{1\dots 2}], [z_{1\dots 4}, \zeta_{1\dots 4}])),$$
$$St(s_k, s_{k + 1}, [z_{12\dots 2}, \zeta_{12\dots 2}], [z_{12\dots 4}, \zeta_{12\dots 4}]) = F_2(St(s_{k - 1}, s_k, [z_{1\dots 2}, \zeta_{1\dots 2}], [z_{1\dots 4}, \zeta_{1\dots 4}])),$$
$$St(s_k, s_{k + 1}, [z_{13\dots 2}, \zeta_{13\dots 2}], [z_{13\dots 4}, \zeta_{13\dots 4}]) = F_3(St(s_{k - 1}, s_k, [z_{1\dots 2}, \zeta_{1\dots 2}], [z_{1\dots 4}, \zeta_{1\dots 4}])), \quad \mbox{and}$$
\begin{equation}
    St(s_k, s_{k + 1}, [z_{14\dots 2}, \zeta_{14\dots 2}], [z_{14\dots 4}, \zeta_{14\dots 4}]) = F_4(St(s_{k - 1}, s_k, [z_{1\dots 2}, \zeta_{1\dots 2}], [z_{1\dots 4}, \zeta_{1\dots 4}]))
\end{equation}
with $p = z_1$ in the definition of transformation $F$, and where $\dots$ denote $k - 2$ ones in a row if the word ends with $2$, and $k - 2$ in a row if the word ends with $4$. The figure $St(s_{k - 1}, s_k, [z_{1\dots 2}, \zeta_{1\dots 2}], [z_{1\dots 4}, \zeta_{1\dots 4}])$ is tiled with $T1$ and $T2$ by the induction hypothesis. So to tile $St(s_k, s_{k + 1}, [z_{11\dots2}, \zeta_{11\dots2}], [z_{14\dots 4}, \zeta_{14\dots 4}])$, we need to check again that 
$$F_1 = \zeta_{11\dots 41}\zeta_{11\dots 4}\zeta_{12\dots 1}\zeta_{12\dots 2}\zeta_{12\dots 21}\zeta_{12\dots 14}\zeta_{12\dots 13}\zeta_{12\dots 12}\zeta_{12\dots 11}\zeta_{11\dots 44}\zeta_{11\dots 43}\zeta_{11\dots 42},$$
where $\dots$ mean $k - 2$ "4"s if the word begins with $11$ or $k - 2$ "1"s if the word begins with $12$, is similar to $T2$, and that 
$$F_3 = \zeta_{12\dots 41}\zeta_{12\dots 4}\zeta_{1311}\zeta_{13\dots 2}\zeta_{13\dots 21}\zeta_{13\dots 14}\zeta_{13113}\zeta_{13\dots 12}\zeta_{13\dots 11}\zeta_{12\dots 44}\zeta_{12\dots 43}\zeta_{12\dots 42},$$
where $\dots$ mean $k - 2$ "4"s if the word begins with $12$ or $k - 2$ "1"s if the word begins with $13$, is similar to $T1$. The figure
$$F_2 = \zeta_{13\dots 41}\zeta_{13\dots 4}\zeta_{14\dots 1}\zeta_{14\dots 2}\zeta_{14\dots 21}\zeta_{14\dots 14}\zeta_{14\dots 13}\zeta_{14\dots 12}\zeta_{14\dots 11}\zeta_{13\dots 44}\zeta_{13\dots 43}\zeta_{13\dots 42}$$
(with "4"s after $13$ and "1"s after $14$) is symmetric to $F_1$, $F_2 = -\overline{F_1},$ so it is also similar to $T2$, is $F_1$ is. One follows the scheme above to check that $F_1$ is similar to the tile $T2$ above it in the stripe between $s_{k - 1}$ and $s_k$, $St(s_{k - 1}, s_k, [z_{11\dots 4}, \zeta_{11\dots 4}], [z_{12\dots 2}, \zeta_{12\dots 2}])$ (where $\dots$ means $k - 3$ "4"s or "1"s in a row) and that the same is true for $F_3$ and the tile $T1$ above it.


\section{Mollified curves and tiles}

The large scale net we constructed cannot yet be a part of the system of curves we reconstruct our coefficient $a$ from, since every green curve from our family $\{s_k\}$ is not smooth at countably many points, and moreover whenever a red curve crosses a green one, they do not meet orthogonally as required (because tangents at these points are not even defined). We tackle this by correcting green curves around each of these bad points. We will do so in a way that preserves the self-similarity, the symmetries, and the length conditions (including property b) from the definition of the curve that joins $\zeta_{114}$ and $\zeta_{122}$, and everything what follows from it).

Consider the green curve $s_2$. Choose a radius $0 < r << 1$ such that the property c) for the curve $(\zeta_{114}, \zeta_{122})$ is satisfied for the radius $r_2 = l^2 r$. Mollify the curve inside the balls with radius $r_2$ centred at points $\zeta_{133}$, $\zeta_{114}$ and $\zeta_{121}$, keeping it symmetric inside those balls with respect to the axis $(\zeta_{1131}, \zeta_{113})$, $(\zeta_{1141}, \zeta_{114})$ and $(\zeta_{1211}, \zeta_{121})$ respectively. Mollify as well the parts of the curve $(\zeta_{114}, \zeta_{121})$ around the points where it is not smooth (if needed, make $r$ smaller, so that $r_2 < \frac{1}{10}\min_{x_i, x_j}{|x_i - x_j|}$, where by $x_i$ we denote the points of non-smoothness). Evidently there is some freedom about how we mollify $s_2$ around each of the enlisted points. The exact procedure does not really matter, but the result has to satisfy the following restriction. It is very important for our construction that the pieces of the curve $(\zeta_{113}, \zeta_{114})$ and $(\zeta_{114}, \zeta_{121})$ have the same length. We want to keep this property for the mollified versions of the curve.

\begin{figure}[h]
    \centering
    \includegraphics[scale = 0.15]{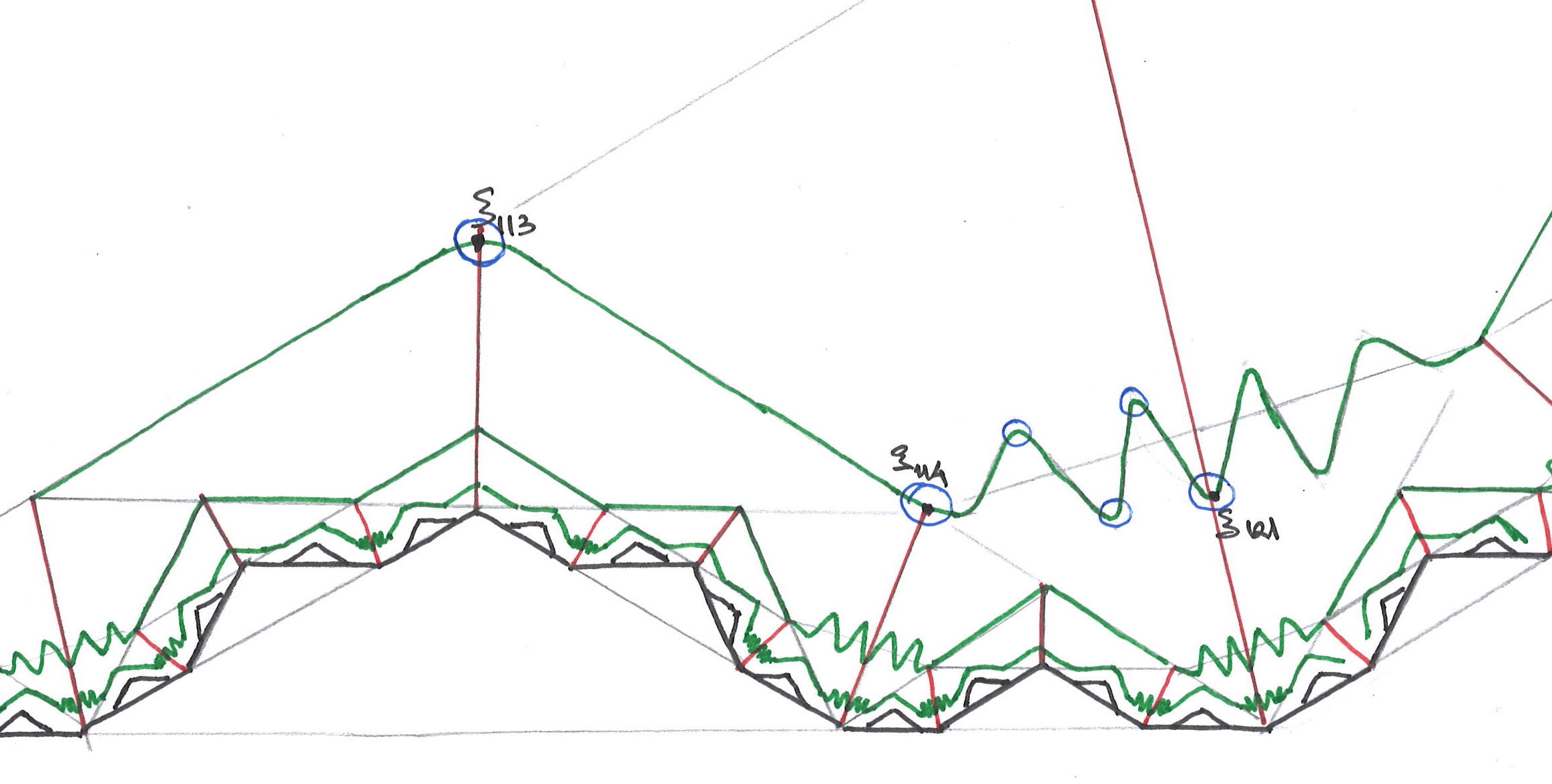}
    \caption{The mollified curve $s_2$ }
    \label{fig:mollified s2}
\end{figure}

Therefore we prefer to write explicitly what we do. Take the ball $B(\zeta_{113}, r_2)$. The curve $s_2 \cap B(\zeta_{113}, r_2)$ is symmetric with respect to $[\zeta_{1131}, \zeta_{113}]$. Consider $s_2 \cap B(\zeta_{113}, r_2)$ as a graph of a symmetric function $f$ with the line $L$ orthogonal to $[\zeta_{1131}, \zeta_{113}]$ as a domain of definition: $f: B(\zeta_{113}, r_2) \cap L \to B(\zeta_{113}, r_2)$. Choose a mollifier $\phi_2: \R \to \R$ such that $\phi_2 = 0$ outside $[-2/3 r_2, 2/3 r_2]$, $\phi_2 = 1$ inside $[-1/3 r_2, 1/3 r_2]$, smooth and takes intermediate values on $[-2/3 r_2, -1/3 r_2]$ and $[1/3 r_2, 2/3 r_2]$. Replace the curve $s_2 \cap B(\zeta_{113}, r_2)$ with the graph of function $f * \phi_2$. Repeat the same for the ball around $\zeta_{114}$. For the ball around $\zeta_{121}$ and for all the non-smoothness points on the part of the curve $(\zeta_{114}, \zeta_{121})$ we do a bit differently. First we precise that, since we do not have the axis of symmetries prescribed to all the other non-smoothness points on $(\zeta_{114}, \zeta_{121})$ apart from $\zeta_{121}$, as a domain of definition of the function $f$, which represents $s_2$ locally, we take the line $(\zeta_{114}, \zeta_{122})$ (recall the property d) when constructing $s_2$ between $\zeta_{114}$ and $\zeta_{122}$). There exists a function $\psi_2: \R \to \R$ such that it is zero outside, say, $[-2/3 r_2, 2/3 r_2]$, one in a small neighbourhood of $0$, takes intermediate values in between, smooth, and, moreover, once we replaced $s_2$ in all the balls around $\zeta_{121}$ and all the non-smoothness points on $(\zeta_{114}, \zeta_{121})$ by $f * \psi_2$, the mollified versions of $(\zeta_{113}, \zeta_{114})$ and $(\zeta_{114}, \zeta_{121})$ have the same length. So this is what we choose as a smoothing procedure for the neighbourhoods of $\zeta_{121}$ and all the non-smoothness points on $(\zeta_{114}, \zeta_{121})$. See figure $\ref{fig:mollified s2}$. 

Once we are done with $\zeta_{113}$, $\zeta_{114}$ and $\zeta_{121}$, we can either define explicitly the mollification at the rest of the non-smoothness points of the family of the green curves, or reconstruct the mollification of rest of $s_2$ and the family $\{s_k\}_{k > 2}$ by the symmetries of the system, thanks to transformation $F$.

The explicit mollification goes like this. Suppose we have a point of non-smoothness $\zeta_w$ marked on $s_k$. Consider $s_k \cap B(\zeta_w, l^k r)$ as a graph of a symmetric function $f$ with the line orthogonal to $[\zeta_w, \zeta_{w1}]$ as a domain of definition. Define $\phi_k(x) = \phi_2(x/l^{k - 2})$ and $\psi_k(x) = \psi_2(x/l^k r)$. Replace the curve $s_k \cap B(\zeta_{w}, l^k r)$ with the graph of function $f * \phi_k$ or $f * \psi_k$, if $w$ ends on $21$ or $41$, or on $11$ and is not of the form $31\dots 1$ (where we have some ones instead of $\dots$). Or, in other words, we replace with $f * \phi_k$, if at least one part of the curve $s_2$ around $\zeta_w$, $(\zeta_{w - 1}, \zeta_w)$ or $(\zeta_w, \zeta_{w + 1})$, is a segment, and by $f * \psi_k$ otherwise.  For the rest of non-smoothness points $p$, we do the same replacement by $f * \psi_k$, with the only difference with the choice of the domain of definition for $f$ in the procedure. By the construction such points $p$ are situated on a curve which joins $\zeta_{w}$ and $\zeta_{w + 2}$. We choose the line $(\zeta_w, \zeta_{w + 2})$ as the domain of definition.

If we are to use the symmetries of the construction, here is what we do. Fist we complete the definition of the mollified $s_2$ curve. We adopt for it the local notation $\widetilde{s_2}$. We reflect the mollification of the neighbourhood of $\zeta_{114}$ with respect to $(\zeta_{113}, \zeta_{1131})$ to get the mollification of the neighbourhood of $\zeta_{112}$. Then we reflect the mollified curve $(\zeta_{112}, \zeta_{113}, \zeta_{114}, \zeta_{121})$ with respect to the axis $(\zeta_{121}, z_{121})$. This is indeed the mollification of the initial curve $s_2$, thanks to the property c) from the definition of $(\zeta_{114}, \zeta_{122})$: $s_2$ in the small neighbourhoods of $\zeta_{122}$ is symmetric with respect to $(\zeta_{123}, \zeta_{1231})$ to $s_2$ in the small neighbourhood of $\zeta_{124}$. Last, to mollify $s_2$ in the neighbourhood of $\zeta_{131}$, we reflect the mollification of the neighbourhood of $\zeta_{123}$ with respect to $(\zeta_{124}, \zeta_{1241})$. Then we reflect the new curve with respect to the axis $\{x = 0\}$.

To define the mollified version of the curve $s_k$, $\widetilde{s_k}$, we apply the inductive definition $\eqref{s_k}$ and $\eqref{complete_s_k}$ to the new curve $\widetilde{s_2}$ instead of the non-smooth version $s_2$. Clearly this gives the same result as the explicit mollification procedure given above, in particular because of the property c) in the definition of $(\zeta_{114}, \zeta_{122})$. From now on we will refer to the new smooth versions of the green curves as $\{s_k\}$: $\{s_k\} := \{\widetilde{s_k}\}$. Of course the smooth version does not differ from the old one outside small neighbourhoods of points $\zeta_w$, so overall the new $s_k$ diverges only slightly from the initial curve. 

\begin{figure}[h]
    \centering
    \includegraphics[scale = 0.15]{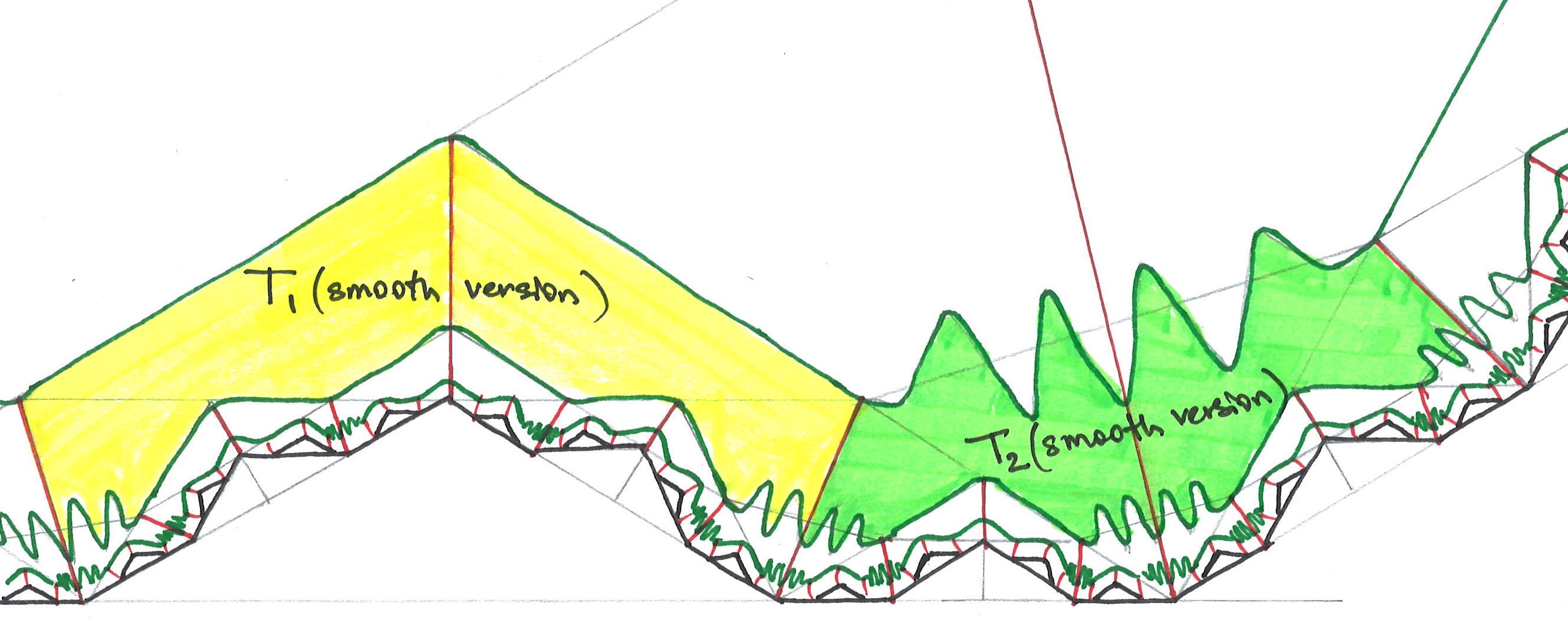}
    \caption{Mollified tiles $T1$ and $T2$}
    \label{fig:mollified tiles}
\end{figure}

But still, consequently, our tiling of the neighbourhood of $S_\alpha|_{[z_1, z_2]}$ is no longer compatible with the large scale net of curves. So we redefine the tiles $T1$ and $T2$ as in the previous section, but whenever we refer in it to a piece of a green curve joining $\zeta_{w1}$ with $\zeta_{w2}$, we refer to the new version of this green curve. This is the same as if we would have applied the explicit mollification procedure to the boundaries of the tiles. Observe that these new smooth figures still tile the space: the proof is word for word the same as in the previous Section. We will now call $T1$ and $T2$ the smooth version of the initial tiles. Note that the new tiles also differ only slightly from the old versions. See figure $\ref{fig:mollified tiles}$

\section{Completing the net and checking the main net property}

We start this section with finally explaining how to build the intermediate red and green curves inside (smooth) tiles $T1$ and $T2$. Each tile has two green and two red ``sides'' of the boundary. We will call the green sides lower side or upper side, and red sides -- left side or right side, according to the layout on the plane of the representatives inside the stripe between $s_2$ and $s_3$. Namely, for the model tile $T1$ we described in the beginning of Section 4, the green sides are $(\zeta_{112}, \zeta_{113}, \zeta_{114})$, upper, and $(\zeta_{1121}, \zeta_{1122}, \zeta_{1123}, \zeta_{1124}, \zeta_{1131}, \zeta_{1132}, \zeta_{1133}, \zeta_{1134}, \zeta_{1141})$, lower, with two red sides which are segments -- left, $[\zeta_{1121}, \zeta_{112}]$, and $[\zeta_{1141}, \zeta_{114}]$, right. The model tile $T2$ from Section 4 has two green sides $(\zeta_{114}, \zeta_{121}, \zeta_{122})$, upper, and $(\zeta_{1141}, \zeta_{1142}, \zeta_{1143}, \zeta_{1144}, \zeta_{1211}, \zeta_{1212}, \zeta_{1213}, \zeta_{1214}, \zeta_{1221})$, lower, and two red -- the left $[\zeta_{1141}, \zeta_{114}]$, and $[\zeta_{1221, \zeta_{122}}]$, the right. We denote by $x_1$ the length of the upper green side of a tile, and by $x_2$ the length of the lower green side of a tile; by the construction it is true that the lower and upper side lengths are the same for the tiles of both types of the same scale (but it is going to be particularly important later). 

We construct the intermediate red curves first. Recall that each tile has an axis of symmetry which divides it into the left half and the right half. Treat a point on the upper side as a $t_1 \in [0, x_1]$, and a point on the lower side as a $t_2 \in [0, x_2]$. We draw smooth curves which enter orthogonally the lower and the upper side and join $t_1$ with the point $t_2 = \frac{x_2}{x_1} t_1$ on the left half of the tile: for $t_1 < x_1/2$. For example, we can draw segments between those two points, and then modify them a little bit near the green sides so they enter them orthogonally (and so that all the curves stay at least $C^4$-smooth). Then we reflect the curves on the left half of the tile to the right half of the tile. Observe that they will still enter lower and upper sides orthogonally, and that they joint a point $t_1$ on the upper side with the point $\frac{x_2}{x_1} t_1$ on the lower one.

\begin{figure}[h]
    \centering
    \includegraphics[scale = 0.2]{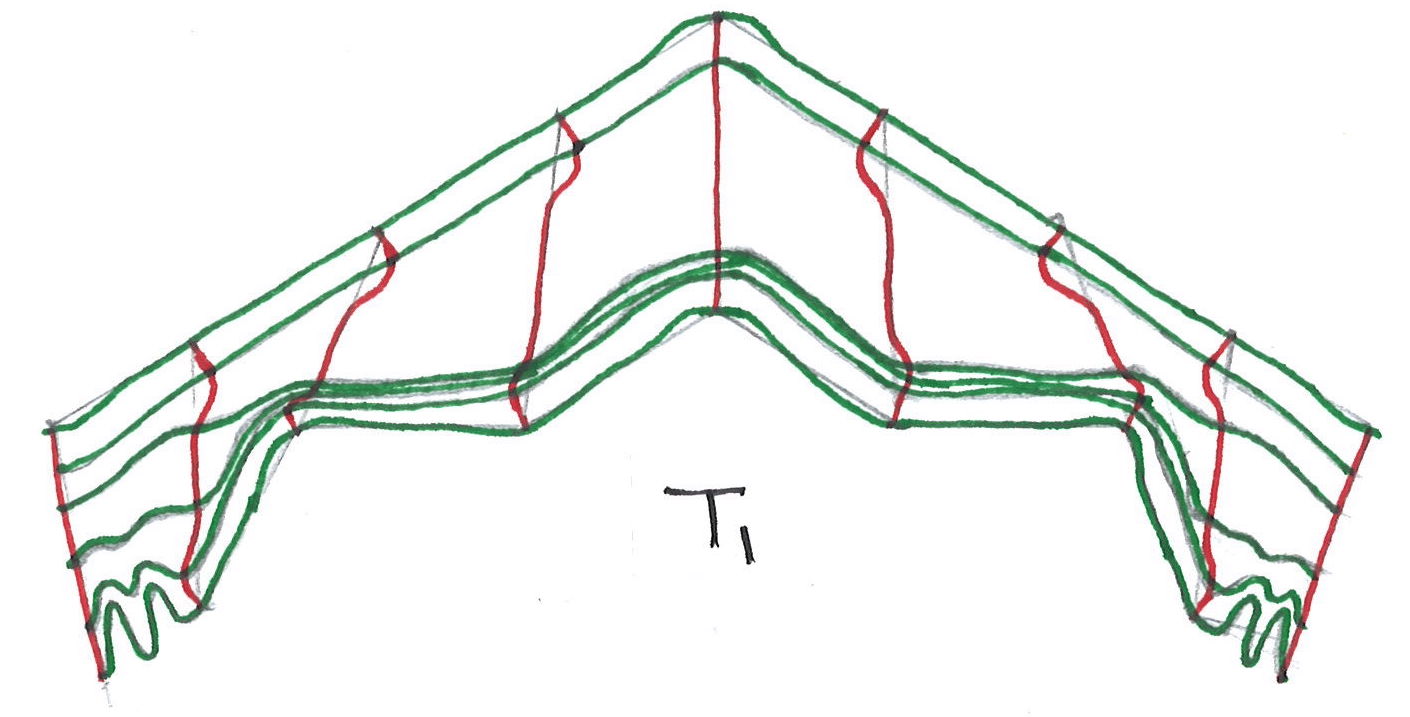}
     \caption{Curves inside $T1$}
    \label{fig:curves T1}
\end{figure}

\begin{figure}[h]
    \centering
     \includegraphics[scale = 0.2]{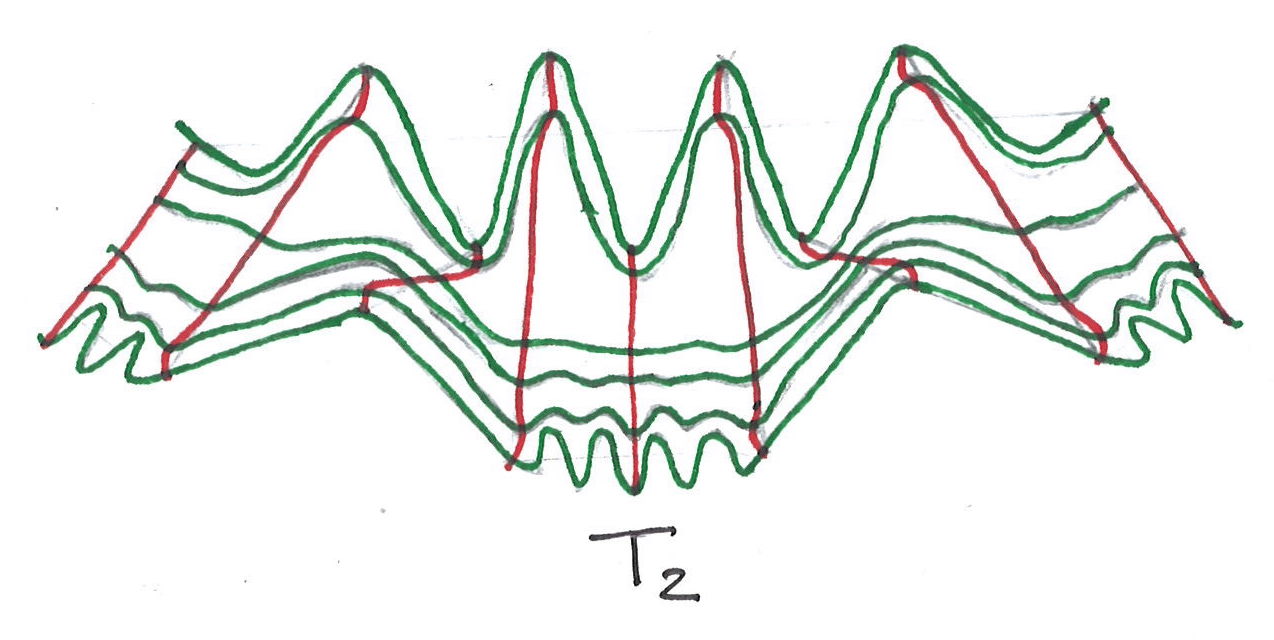}
    \caption{Curves inside $T2$}
    \label{fig:curves T2}
\end{figure}

By the procedure described in Section 2, we reconstruct the green curves in the tile (the procedure reconstructs red curves from the green, but of course the same can be done the other way round, and if a function $v$ is conjugate to $u$, then $u$ is conjugate to $v$). Observe that all the green curves are symmetric with respect to the axis of symmetry of the tile, thanks to the construction of the red curves. See figures $\ref{fig:curves T1}$ and $\ref{fig:curves T2}$.

%
%
%
   

Thus the net of red and green curves is built on the whole triangular neighbourhood of $S_\alpha$ above $[z_1, z_2]$, apart from the areas near the legs $[z_1, \zeta_{13}]$ and $[\zeta_{13}, z_2]$ (they are composed of all the pieces of type $\zeta_{1\dots 14}\zeta_{1\dots 13}\zeta_{\dots 12}\zeta_{\dots 21}$ and $-\overline{\zeta_{1\dots 14}\zeta_{1\dots 13}\zeta_{\dots 12}\zeta_{\dots 21}}$ with ones instead of $\dots$): we equip all the tiles of the same type with the net of curves as above. Recall that, if figures $F1$ and $F2$ are tiles of the same type, by construction we have $F2 = g(F1)$ where $g = h \circ m$ with a homothety $h$ and a motion $m$, and we use the map $g$ to transfer the coloured net to $F2$.

~\ 

It is time we check the main property $\eqref{mainproperty}$ for the net we built above $S_\alpha$ constructed on $[z_1, z_2]$. First, we need to choose a green curve and parameterise it. Choose the piece of $s_2$ we constructed and parameterise it naturally (with length). This restricts for now the neighbourhood of our considerations to the tiled region above $S_\alpha |_{[z_{112}, z_{144}]}$, since $s_2$ ends on the left with $\zeta_{112}$ and on the right with $\zeta_{144}$, and the red curves which join the defined piece of $s_2$ with $S_\alpha$ end on $z_{112}$ and $z_{144}$.

Second, we need to complete the labelling. We have just labelled our red curves, it is now the turn of the green ones. We start with declaring the labelling for the discrete family $\{s_k\}$. We put
\begin{equation}\label{green_label}
    s_k = l^{kd},
\end{equation}
where $d$ is the Hausdorff dimension of the snowflake $S_\alpha$. We label the intermediate curves $s$ the following way. Assume that a green curve $s$ crosses the left lid of a tile in a stripe between $s_k$ and $s_{k + 1}$. Identify the left lid of the tile with a segment $[0, h_k]$; one can compute $h_k$ explicitly, but we do not need to do it. The curve $s$ crosses the lid somewhere at a point $y \in [0, h_k]$. Put 
\begin{equation}\label{green_label2}
    s = l^{kd}\frac{y}{h_k} + l^{(k + 1)d}\frac{h_k - y}{h_k}.
\end{equation}
This definition is coherent, meaning, it does not depend on the tile we chose, since all the red sides of tiles in the stripe $St(s_k, s_{k + 1})$ have the same length, and if a green curve starts from a point $y$ on the left lid of a tile, it arrives at the symmetric point $y$ on the right lid of the tile. So we can identify all of them with $[0, h_k]$, and $s$ crosses each of them at the same point $y \in [0, h_k]$. 

Labelling done, we first check property $\eqref{mainproperty}$ for the tiles that are adjacent to $s_2$, or, in other words, fill in the stripe $St(s_2, s_3, [z_{112}, \zeta_{112}], [z_{144}, \zeta_{144}])$. It suffices to consider a model tile $T1$ $$\zeta_{112}\zeta_{113}\zeta_{114}\zeta_{1141}\zeta_{1134}\zeta_{1133}\zeta_{1132}\zeta_{1131}\zeta_{1124}\zeta_{1123}\zeta_{1122}\zeta_{1121},$$
and a model tile $T2$ $$\zeta_{114}\zeta_{121}\zeta_{122}\zeta_{1221}\zeta_{1214}\zeta_{1213}\zeta_{1212}\zeta_{1211}\zeta_{1144}\zeta_{1143}\zeta_{1142}\zeta_{1141}.$$ 
Without loss of generality, concentrate on the first tile. Choose a cell of the net containing $z \in T1$ which consist of pieces of two green curves labelled $s_a$ and $s_b$, and two red segments labelled, say, $\theta_1$ and $\theta_2$. Denote $w_G(z)$ the vector of speed on a green curve at a point $z$, and $w_R(z)$ the vector of speed on a red curve at a point $z$. We know that there are constants $c_1$ and $c_2$ such that $c_1^{-1} \le w_G(z) \le c_1$ and $c_2^{-1} \le w_R(z) \le c_2$ everywhere on $T1$. Plus, $s_a$ and $s_b$ can be represented as $s_a(t) = \ilim_0^t{w_{G_a}(z(t))dt}$, $s_b(t) = \ilim_0^t{w_{G_b}(z(t))dt}$. Thanks to the nice dependence of the solution on the initial data, the distance between the curves $s_a$ and $s_b$ is comparable to $|s_a - s_b|$. The same is true about $\theta_1$ and $\theta_2$. So we have 
$\frac{\dist_G}{\dist_R} \approx \frac{|s_a - s_b|}{|\theta_1 - \theta_2|}$, therefore property $\eqref{mainproperty}$ follows.

Now pick a tile $T$ from the stripe $St(s_k, s_{k + 1})$, and a point $z$ in it. Pick a cell containing $z$ composed by two green curves $a$ and $b$, $a, b \in [0, h_k]$, and two red curves $\Tilde{\theta_1}$ and $\Tilde{\theta_2}$, $\Tilde{\theta_1}, \Tilde{\theta_2} \in [0, x_k]$, where $[0, x_k]$ parameterises the upper lid of the tile. Note that $a, b, \Tilde{\theta_1}$ and $\Tilde{\theta_2}$ is a ``local'' labelling, and our goal here is to reconstruct the global one $s_a, s_b, \theta_1$ and $\theta_2$. From the previous paragraph we know that the ratio between the ``green'' distance $\dist_G$ and the ``red'' distance $\dist_R$ is approximately $\frac{|a - b|}{|\Tilde{\theta_1} - \Tilde{\theta_2}|}.$ Therefore 
\begin{equation}\label{loc_mainprop}
  \frac{\dist_G |\theta_1 - \theta_2|}{\dist_R |s_a - s_b|} \approx \frac{|a - b||\theta_1 - \theta_2|}{|s_a - s_b| |\Tilde{\theta_1} - \Tilde{\theta_2}|}.  
\end{equation}
Denote labels of the lids of the tile $\theta_l, \theta_r$, $s_u$ and $s_d$, the indices stand for left, right, up and down, $s_u = s_k$ and $s_d = s_{k + 1}$. By the construction of intermediate curves we have
$$|s_a - s_b| = |s_{k + 1} - s_k|\frac{|a - b|}{h_k} \quad \mbox{and} \quad |\theta_1 - \theta_2| = |\theta_l - \theta_r|\frac{|\Tilde{\theta_1} - \Tilde{\theta_2}|}{x_k}.$$
So $\eqref{loc_mainprop}$ transforms into
$$\frac{|a - b||\theta_1 - \theta_2|}{|s_a - s_b| |\Tilde{\theta_1} - \Tilde{\theta_2}|} = \frac{h_k |a - b||\theta_l - \theta_r||\Tilde{\theta_1} - \Tilde{\theta_2}|}{x_k |s_{k + 1} - s_k| |a - b| |\Tilde{\theta_1} - \Tilde{\theta_2}|} = \frac{h_k}{x_k} \frac{|\theta_l - \theta_r|}{|s_{k + 1} - s_k|}.$$
The ratio $\frac{h_k}{x_k}$ by the construction is constant for the tiles of the same type, so it is a global constant (which depends on $\alpha$) and does not depend on $k$. Therefore it is left for us to check that
\begin{equation}\label{comp_labels}
    |s_{k + 1} - s_k| \approx |\theta_l - \theta_r|.
\end{equation}
By $\eqref{green_label}$ we have
$$|s_{k + 1} - s_k| = |l^{(k + 1)d} - l^{kd}| = l^{kd}(1 - l^d).$$
Computing $|\theta_l - \theta_r|$ is a little bit more subtle. To do this we need to trace back the history of the red sides of our tile $T$ and find where they cross the curve $s_2$, since its natural parameterization serves as the global labeling of red curves. Luckily we are not interested in the exact values of $\theta_r$ and $\theta_l$, only their difference is important, so the matters are slightly easier than that. By the construction of tiles in the stripe $St(s_k, s_{k + 1})$, the upper green lid has length $2\frac{l^k}{\cos{\alpha}}$, and the lower green lid has length $8\frac{l^{k + 1}}{\cos{\alpha}}$. So the piece of $s_k$ between $\Tilde{\theta_l}$ and $\Tilde{\theta_r}$ in ``local'' coordinates, or $\theta_l$ and $\theta_r$ in the global coordinates, has length $2\frac{l_1^k}{\cos{\alpha}}$. By definition of intermediate red curves, the piece of $s_{k - 1}$ between $\theta_l$ and $\theta_r$ has length $ 2\frac{l_1^k}{\cos{\alpha}} \frac{2\frac{l^k}{\cos{\alpha}}}{8\frac{l^{k + 1}}{\cos{\alpha}}} = 2\frac{l^k}{\cos{\alpha}} (4l)^{-1}$. Tracing back by induction the length of pieces of curves between $\theta_l$ and $\theta_r$, we'll have that on $s_2$ the length between $\theta_l$ and $\theta_r$ is
$$2\frac{l^k}{\cos{\alpha}} (4l)^{-(k - 2)} \approx 4^{-k}.$$
Recall that we have $|s_{k + 1} - s_k| \approx l^{kd}$. The Hausdorff dimension $d$ of the snowflake is $\frac{\ln{(4)}}{\ln{(1/l)}}$, so $l^{-d} = 4$, which implies that $\eqref{comp_labels}$ holds.

\section{The proof of the main theorem}

Before we finally prove Theorem $\ref{mainthm1}$, we need to complete the net of green and red curves on the whole upper-space for both types of our snowflake $S_\alpha$ announced in the introduction -- $S_\alpha^1$ and $S_\alpha^2$. 

We start with $S_\alpha^1$, because, after we deal with it, it is going to be clear how to deal with $S_\alpha^2$ as well. Recall the definition of $S_\alpha^1$ from p. 5. First, take as the initial set $I$ the snowflake $S_\alpha|_{[-1/2, 1/2]}$, put $(S_\alpha)_1 = I = S_\alpha|_{[-1/2, 1/2]}$. Define $(S_\alpha)_2 = H_l((S_\alpha)_1)$, where $H_l$ is the scaling transform with the coefficient $1/l$ defined by $\eqref{l}$ and the center $-1/2$. By induction, define $(S_\alpha)_n = H_l((S_\alpha)_{n - 1})$. Finally, we define
$$S_\alpha^1|_{\{x \ge -1/2\}} = \cup_{n \ge 1}{(S_\alpha)_n} \quad \mbox{and}$$
\begin{equation}\label{S_alpha_one}
    S_\alpha^1 = S_\alpha^1 |_{\{x \ge -1/2\}} \cup \left[-\overline{\left(S_\alpha^1 |_{\{x \ge -1/2\}} + \frac{1}{2}\right)} - \frac{1}{2}\right]
\end{equation}
See figure $\ref{fig:S 1 alpha}$.

\begin{figure}[h]
    \centering
    \includegraphics[scale= 0.09]{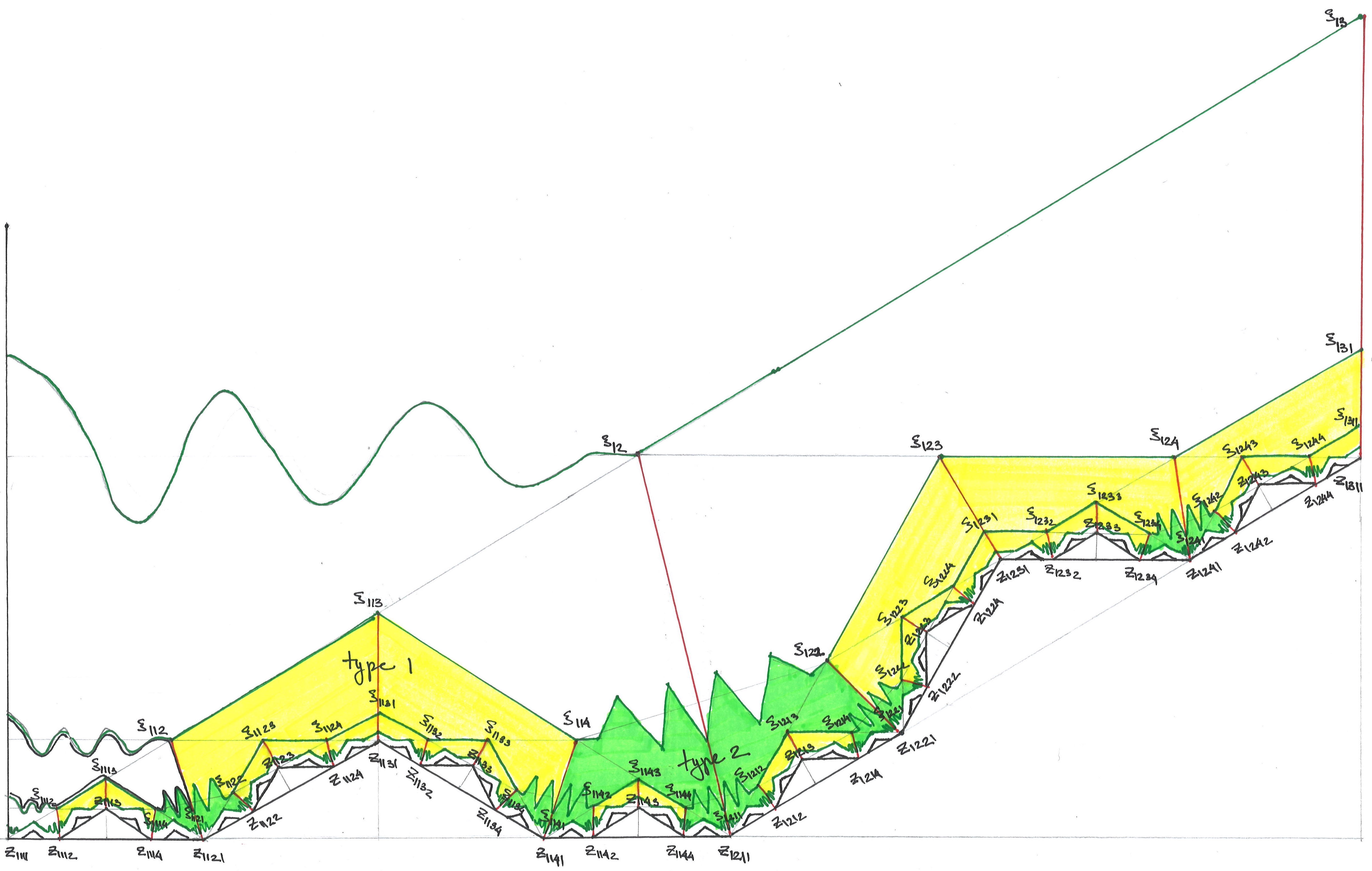}
    \caption{Globally defined green curves near the segment $[-1/2, 0]$}
    \label{fig:global green curves}
\end{figure}

This definition also extends our green and red curves $\{s_k\}$ and $\{\theta_w\}$. Indeed, for example, the set $\cup_{n \ge 0}H_l^n(s_{k + n})$ extends the green curve $s_k$ to the whole half-plane $\{x \ge 0\}$ (to be more precise, to the right of the point $\zeta_{1 \dots 12}$ with $k - 2$ "1"s). To extend the curve $s_k$ to the left of the interval $[-1/2, 1/2]$, we just reflect the whole picture with respect to $\{x = -1/2\}$, as in $\eqref{S_alpha_one}$. If a point $p$ lies on the half-space $\{x \ge -1/2\}$, we refer to the point symmetric to it with respect to $\{x = -1/2\}$ as to $s(p)$. Note that we still need to complete the definition of green curves $s_k$ by joining $\zeta_{1\dots 12}$, where $\dots$ are $k - 2$ "1"s, with $s(\zeta_{1\dots 12})$. We do something similar to when we were inventing the curve $(\zeta_{114}, \zeta_{121})$. Draw a curve joining $\zeta_{12}$ with $s(\zeta_{12})$ such that
\begin{enumerate}
    \item[a)] it is symmetric with respect to $\{x = -1/2\}$ and at $x = -1/2$ it has horizontal tangent,
    \item[b)] smooth,
    \item[c)] has twice the length of $[\zeta_{12}, \zeta_{13}]$, $2\frac{l}{\cos{\alpha}}$,
    \item[d)] in small (radius much smaller than $l^2$) ball centred at $\zeta_{12}$, the curve is symmetric to the segment $[\zeta_{12}, \zeta_{13}]$ with respect to the axis $(z_{12}, \zeta_{121})$,
    \item[e)] the curve stays inside a small neighbourhood of the line $y = \Im(\zeta_{12})$, where small means, say, at distance less than $\frac{\Im(\zeta_{12}) - \Im(\zeta_{113})}{2}$ (so that it does not intersect the curve $s_2$).
\end{enumerate}
To join $\zeta_{1\dots 12}$ with $s(\zeta_{1\dots 12})$ on the curve $s_k$, put
\begin{equation}\label{type3_curve}
    (s(\zeta_{1\dots 12}), \zeta_{1\dots 12}) = (s(\zeta_{112}, \zeta_{112}) - z_1)l^{k - 2} + z_1.
\end{equation}
See figure $\ref{fig:global green curves}$ for the sketch of the curves to the right of the point $z_1$. 

\begin{figure}[h]
    \centering
    \includegraphics[scale= 0.09]{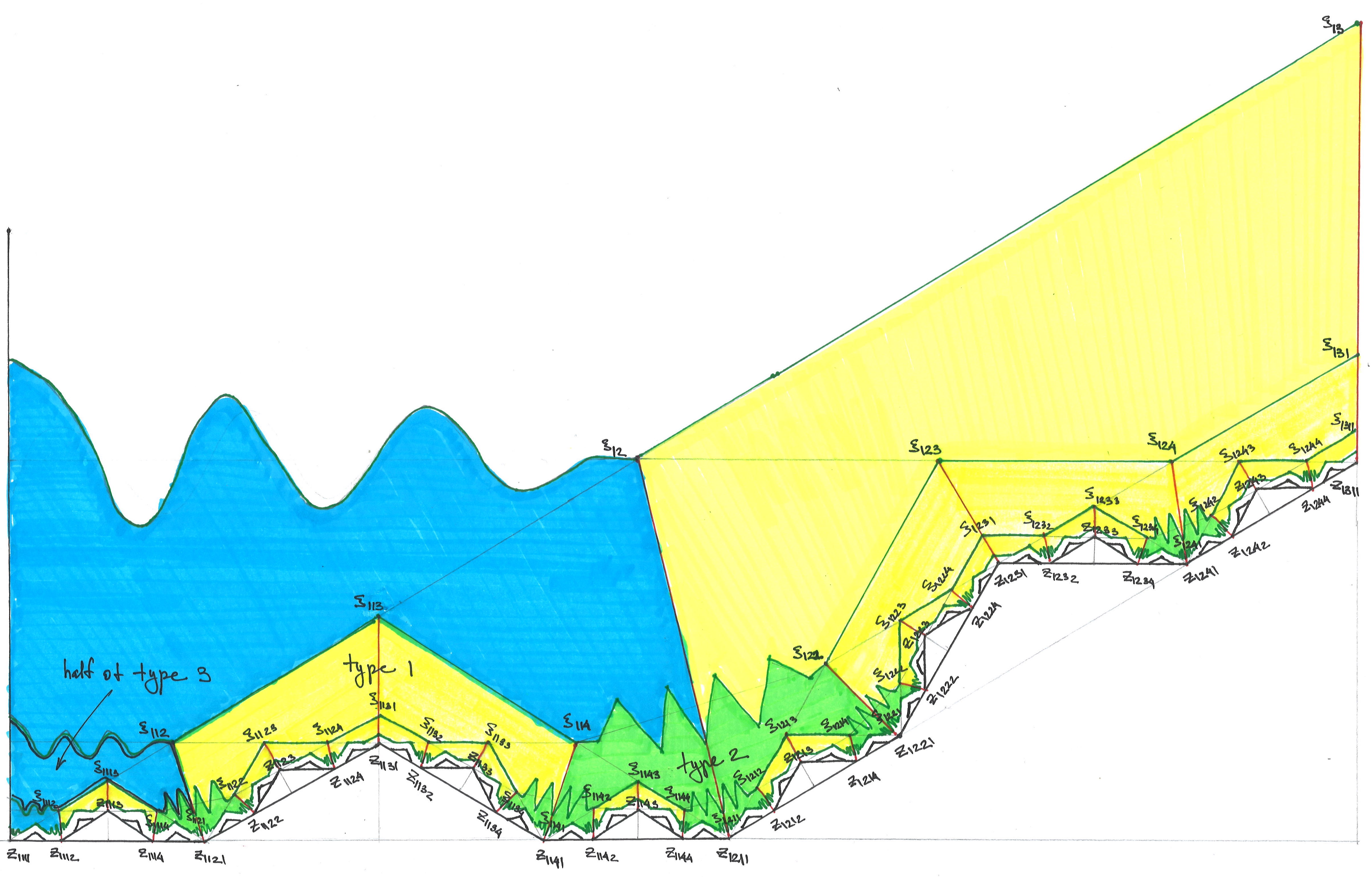}
    \caption{Global tiling near the segment $[-1/2, 0]$, tiles $T3$}
    \label{fig:global tiling}
\end{figure}

Here is how the upper-half space above $S_\alpha^1$ is tiled. First, we use the scaling $H_l$ to extend the tiling to the right of the interval $[-1/2, 1/2]$. Denote by $U_1 = U$ the neighbourhood of $S_\alpha|_{[-1/2, 1/2]}$ we tiled and filled with a net of curves in Sections 3-6. By the inductive definition $\eqref{S_alpha_one}$, on the first step we tile $U_2 = H_l(U_1)$ by images of tiles under $H_l$. On step $n$, we tile $U_n = H_l(U_{n - 1})$. To get the tiling to the left of the interval $[-1/2, 1/2]$, we reflect the whole construction above $S_\alpha^1 |_{\{x \ge -1/2 \}}$ with respect to $\{x = -1/2\}$, as in $\eqref{S_alpha_one}$. Thus, the domain $(\R^2 \setminus S_\alpha^1)^{+}$ is tiled except for the sector between the rays $\{y = -\tan{(\alpha)}(x + 1/2) \}$ and $\{y = \tan{(\alpha)}(x + 1/2) \}$, plus the images of $\zeta_{1113}\zeta_{112}\zeta_{1121}\zeta_{1114}$ under $H_l^n, n \ge 0$, and the symmetric with respect to the axis $\{x = -1/2\}$ regions. It is left for us to tile these parts. Denote $\zeta_{1 \dots 1} = \{x = -1/2\} \cap (s(\zeta_{1\dots 12}), \zeta_{1\dots 12})$. We call tiles of type 3 all the regions which are similar to the figure
\begin{equation}
    T3 = s(\zeta_{112})\zeta_{111}\zeta_{112}\zeta_{1121}\zeta_{1114}\zeta_{1113}\zeta_{1112}\zeta_{1111}s(\zeta_{1112})s(\zeta_{1113})s(\zeta_{1114})s(\zeta_{1121}).
\end{equation}
The whole region below $s_1$ near $z_1 = -1/2$, which was not tiled yet, tiles with the images $H_l^n(T 3), n \in \Z$. We build the net of intermediate red and green curves in the same way we did it inside $T1$ and $T2$. The proof that this piece of net also satisfies $\eqref{mainproperty}$ is analogous to what was done in Section 6. See on figure $\ref{fig:global tiling}$ how the complete tiling looks like between the segment $[-1/2, 0]$ and the curve $s_1$. To get the idea of how the global tiling looks like, one can zoom in around $z_1$ the picture $\ref{fig:global tiling}$.

~\ 

\begin{figure}[h]
    \centering
    \includegraphics[scale=0.2]{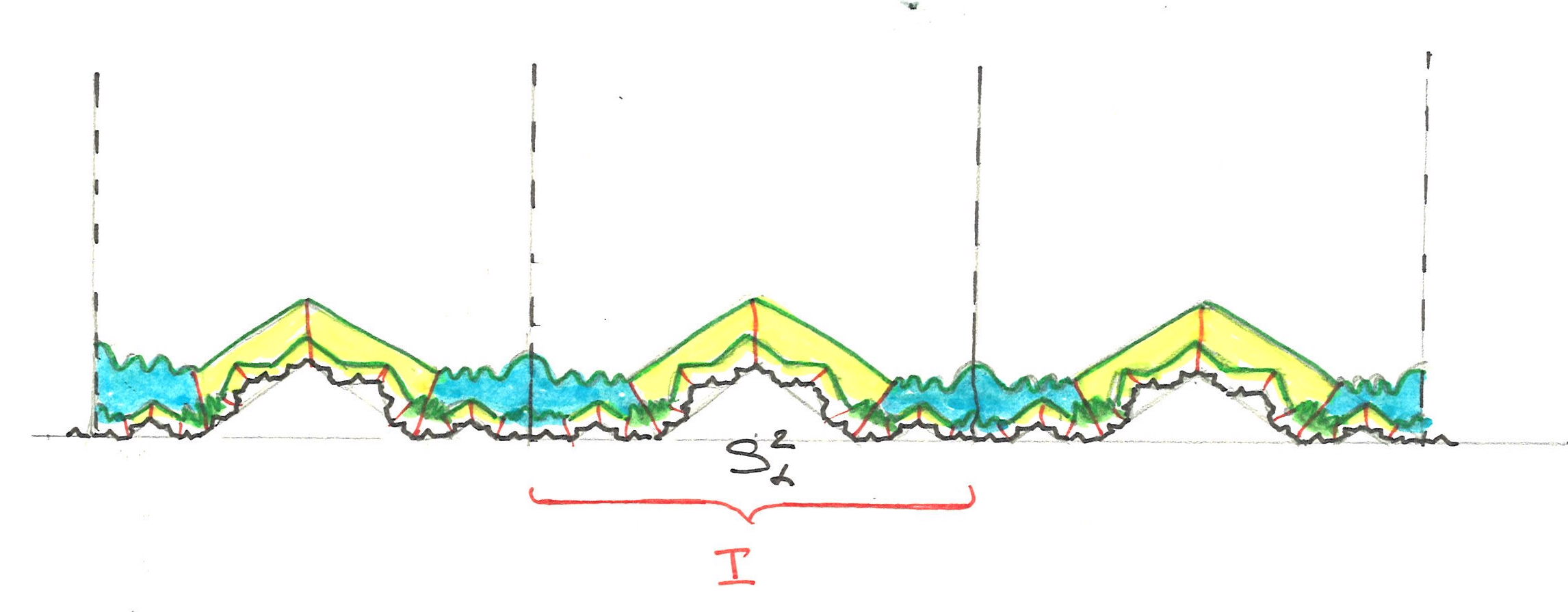}
    \caption{Global tiling above $S_\alpha^2$}
    \label{fig:s2 global tiling}
\end{figure}

We now treat $S_\alpha^2$. Recall that by definition $S_\alpha^2 = \cup_{n \in \Z}(S_\alpha |_{[-1/2, 1/2]} + n)$. We tile the triangular neighbourhood $U_n = U + n$ above $S_\alpha|_{[-1/2, 1/2]} + n$ with figures which tile the neighbourhood $U$ above $S_\alpha|_{[-1/2, 1/2]}$ translated by $n$. Areas above the points with half-integer coordinated are not tiled yet, and the green curves $s_k$ are not defined there as well. We complete the definition of $s_k$ in the same way we did for $S_\alpha^1$ above $z_1$. Again, we use the notation $s(p)$ -- the point symmetric to $p$ with respect to $\{x = -1/2\}$. We join $\zeta_{1\dots 12}$ with $s(\zeta_{1\dots 12})$ as we do for $S_\alpha^1$, see $\eqref{type3_curve}$. For the rest of half-integer points, put
$$s_k|_{((s(\zeta_{1\dots 12}), \zeta_{1\dots 12})) + n} = ((s(\zeta_{1\dots 12}), \zeta_{1\dots 12})) + n.$$
All the areas below the curve $s_1$ around the half-integer points which are not tiled yet can be tiled with $T3$. The intermediate green and red curves are build as above, and $\eqref{mainproperty}$ holds for the resulting net by the same argument again as in Section 6. See figure $\ref{fig:s2 global tiling}$ for a sketch of the global tiling we did. Note that the part of the domain $(\R^2 \setminus S_\alpha^2)^{+}$ above the curve $s_1$ is not tiled yet.

~\ 

Therefore, for $S_\alpha^2$ the coefficient $a$ above the curve $s_1$ is still not defined, so we do it now. Draw a straight horizontal line $s_0$ above $s_1$: one can take $s_0 = \{y = 1\}$. Define $a = 1$ on $\{y \ge 1\}$. Then, chop $s_0$ into intervals of length one, $s_0 = \cup_{n \in \Z}{[-1/2 + n, 1/2 + n]}$. Build the net of green and red curves on stripe between $s_0$ and $s_1$ in the following way. We concentrate first on $St(s_0, s_1) \cap \{-1/2 \le x \le 1/2\}$. Let $t_1 \in [0, 1]$ be a point on $s_0 \cap \{-1/2 \le x \le 1/2\}$, and $t_2 \in [0, 4\frac{l}{\cos{\alpha}}]$ be a point on $s_1 \cap \{-1/2 \le x \le 1/2\}$. Draw red curves which join $t_1$ with $t_2 = 4\frac{l}{\cos{\alpha}} t_1$ and enter $s_0$ and $s_1$ orthogonally. We can do it the same way as when constructing the net inside our tiles $T1, T2$ and $T3$: join $t_1$ and $t_2$ with a straight line, and then correct it a little bit around $s_0$ and $s_1$. Then reconstruct the green curves from the red ones. Note that this net is symmetric with respect to $\{x = 0\}$. Extend the net on the whole stripe $St(s_0, s_1)$ using the definition of $S^2_\alpha$: to get the net inside $St(s_0, s_1) \cap \{-1/2 + n \le x \le 1/2 + n\}$ for $n > 0$, for example, reflect with respect to $\{x = -1/2 + n\}$ the region $St(s_0, s_1) \cap \{-1/2 + (n - 1) \le x \le 1/2 + (n - 1)\}$.



Thus we have the coefficient $a$ continuous and defined everywhere on $\left(\R^2 \setminus S_\alpha\right)^+$. The bounds $\eqref{scalar_coeff}$ still holds for the globally defined $a$ with the same constant $C$.

~\ 

We are ready to return to Theorem $\ref{mainthm1}$, which we state below for the reader's convenience. Let $\mu$ be the Hausdorff measure $\Hd^d|_{S_\alpha}$ on the snowflake $S_\alpha^1$ or $S_\alpha^2$, where $d = \frac{\ln{(4)}}{\ln{(2(1 + \cos{\alpha}))}}$ is the snowflake's dimension. Let $w_L^z$ be the elliptic measure of the operator $L = -\div \; a \nabla$ on the upper half-plane above the snowflake $S_\alpha$ ($S_\alpha^1$ or $S_\alpha^2$) with the coefficient $a$ constructed above. We denote by $w_L^\infty$ the (weak) limit of $w_L^z$ as $z$ tends to infinity. 

\begin{thm}
    There exists an absolute constant (probably dependent on $\alpha$) $C \ge 1$ and a continuous function $a: \left(\R^2 \setminus S_\alpha\right)^+ \to (0, +\infty)$ such that $\eqref{mainproperty}$ holds, and if $w_L^\infty$ denotes the elliptic measure with pole at $\infty$, associated to the operator $L = -\div \; a \nabla$ on the domain $\Omega = \left(\R^2 \setminus S_\alpha\right)^+$, then $w_L^\infty = \mu$. Also,
    \begin{equation}\label{meas_sets}
        C^{-1}\mu(A) \le w_L^z(A) \le C\mu(A)
    \end{equation}
    for all $z$ such that $\delta(x) = \dist(x, S_\alpha) \ge 1, \; A \subset S_\alpha$ measurable.
\end{thm}

\begin{proof}
    Let $a$ be the coefficient constructed above and $G$ the $L$-harmonic function which level lines are the green curves also constructed above. We claim that $G$ is a constant multiple of the Green function of the operator $L$ with pole at infinity: it is $L$-harmonic, positive, vanishes at the boundary and regular enough. The latter is true by our construction: combining the labelling $\eqref{green_label}$ and the fact that $\dist(S_\alpha, s_k) \approx l^k$, we have that
    $$G(z) \sim \dist(z, S_\alpha)^d.$$
    Since the boundary $S_\alpha$ of our domain $\Omega$ is irregular, we will approximate it with already constructed curves $s_k$, so we can work with $G$ more easily. Denote $\Omega_k = \left(\R^2 \setminus s_k\right)^+$, and concentrate on this approximating domain. $G$ is still $L$-harmonic in $\Omega_k$, and therefore $G_k = G - l^{kd}$ is a (constant multiple) of the Green function on $\Omega_k$ with pole at $\infty$. The boundary $s_k$ is smooth, so for the elliptic measure $w_k$ with pole at infinity associated to $L$ on $\Omega_k$ we can write $d w_k = - a g_k d\mu_k$, where $\mu_k$ is the (normalised) arch-length on $s_k$, $g_k = \frac{\partial G_k}{\partial n}$, and $a$ is our scalar coefficient (see \cite{K}, p.6). 

    We claim that the product $ - a g_k$ does not depend on $k$ and is equal to one everywhere on $s_k$. Indeed, suppose we know that $ - a g_2 = 1$ on $s_2$. Let $p$ be a point on the curve $s_k$, $k \neq 2$. Without loss of generality, we can assume that it does not coincide with any of the points $\zeta_w$. Then, $p$ lies on the upper green side of a tile $F$ of the type $T1, T2$ or $T3$. By definition this means that there exists an affine transformation $f = h \circ m$, where $h$ is a scaling with the coefficient $l^{k - 2}$ and $m$ a motion, such that $F = f(T_i)$ with $i = 1, 2$ or $3$. Denote $p' = f^{-1}(p)$. Then $g_k(p) = g_2(p')$, since $G_k(z) + \const = G(f^{-1}(z))$, and the normal vector $n$ at two points coinsides up to the correct scaling. We also have $a(p) = a(p')$ (because locally pictures are the same, but scaled). 
    

    Concerning the product $- a g_2 = - a \frac{\partial G_2}{\partial n}$ on the curve $s_2$, by definition of $a$ we have
    $$- a \frac{\partial G_2}{\partial n} = - \frac{|\nabla v|}{|\nabla G|} \frac{\partial G_2}{\partial n} = \frac{|\nabla v|}{|\nabla G|} |\nabla G| = |\nabla v|,$$
    where $v$ is the conjugated function to $G$, because the curves $s_k$ are level lines of $G$. The vector $\nabla v$ is tangent to the curve $s_2$ by construction. It is left to recall that our global labelling was chosen in such a way that the curve $s_2$ is naturally paratemetised, so $|\nabla v| = 1$. 
    

    This gives us the statement we want, because $\mu_k$ tends (weakly) to $\mu$ (sets $s_k$ tend to $S^\alpha$, meaning, the Hausdorff distance between $s_k$ and $S^\alpha$ tends to zero), and $w_k$ tends (weakly) to $w_L^\infty$. So we have $w_L^\infty = \lim w_k = \lim \mu_k = \mu$.
    
    The double inequality $\eqref{meas_sets}$ follows from the comparison principle.   
\end{proof}

\section{Variants and open questions}

\begin{figure}[h]
    \centering
    \includegraphics[scale = 0.2]{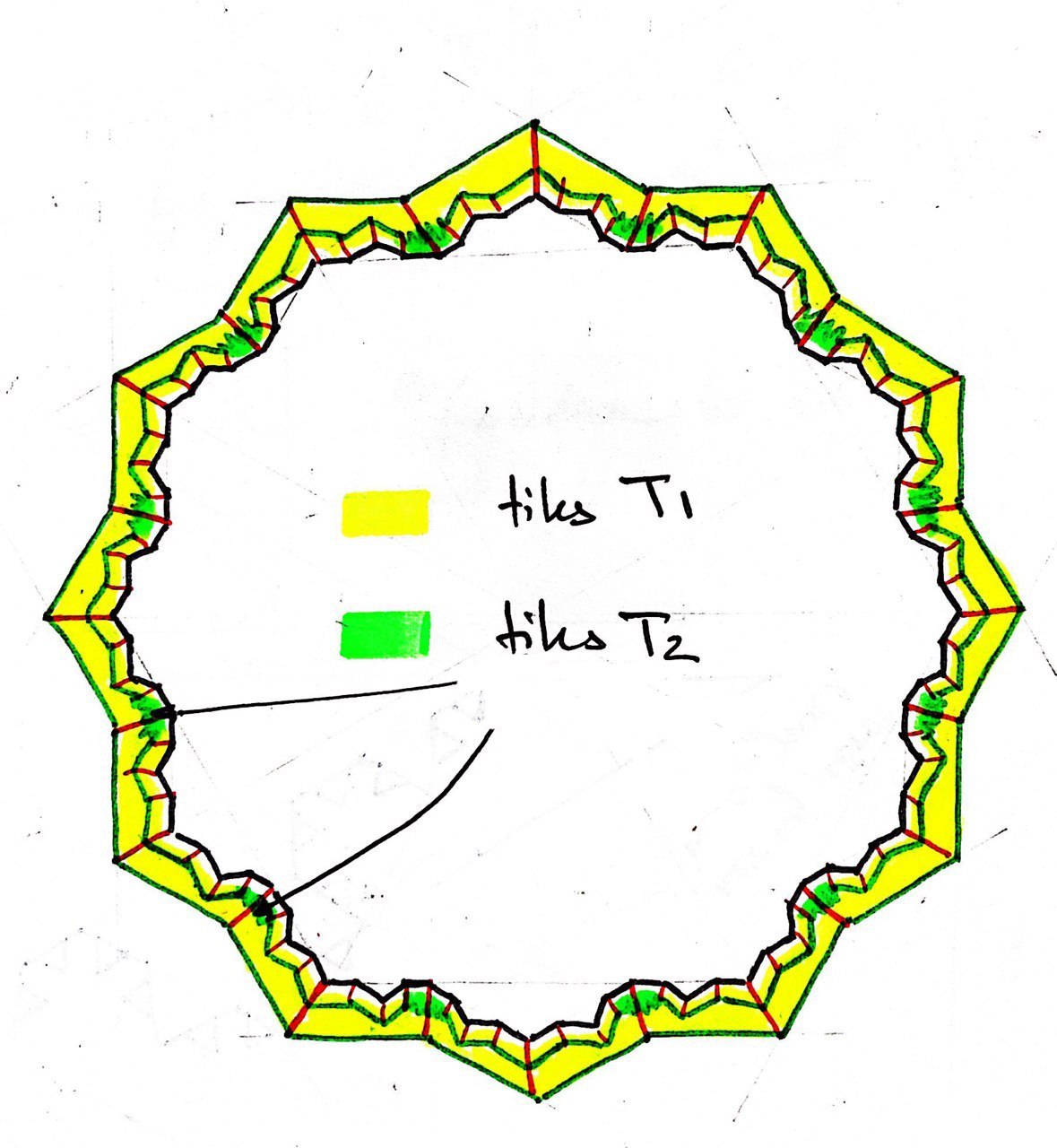}
    \caption{Tiling for a compact snowflake with $\alpha = \frac{\pi}{N}$}
    \label{fig:comp regular tiling}
\end{figure}

Theorem $\ref{mainthm1}$ also covers compact versions of our snowflakes. To construct a first compact version, we could start from, like one does with the classical Koch snowflake, with a regular polygon of angle $\pi - 2\alpha$. If in this case our angle $\alpha$ in the snowflake construction is of the form $\frac{\pi}{N}$, where $N$ is the number of vertices of our regular polygon, we can preserve the symmetries present in the classical construction. The upside is that these symmetries will allow us to eliminate $T3$ from the tiling. One could say that this is more elegant than the unbounded fractal version we have in the previous section, but our original intention was to snow some counterexample on a non-compact unbounded set. Plus, we do not want to limit ourselves to the angles of the form $\frac{\pi}{N}$ only. We do not want to enter the details of how our tiling works in this case, but locally everything works the same as in Sections 3-6, and the extension to the rest of the domain $(\R^2 \setminus S_\alpha)^{+}$ (the non-compact of the two possible options) is quite straightforward. See figure $\ref{fig:comp regular tiling}$ for a sketch of tiling around the compact snowflake $S_\alpha$.


\begin{figure}[h]
    \centering
    \includegraphics[scale = 0.2]{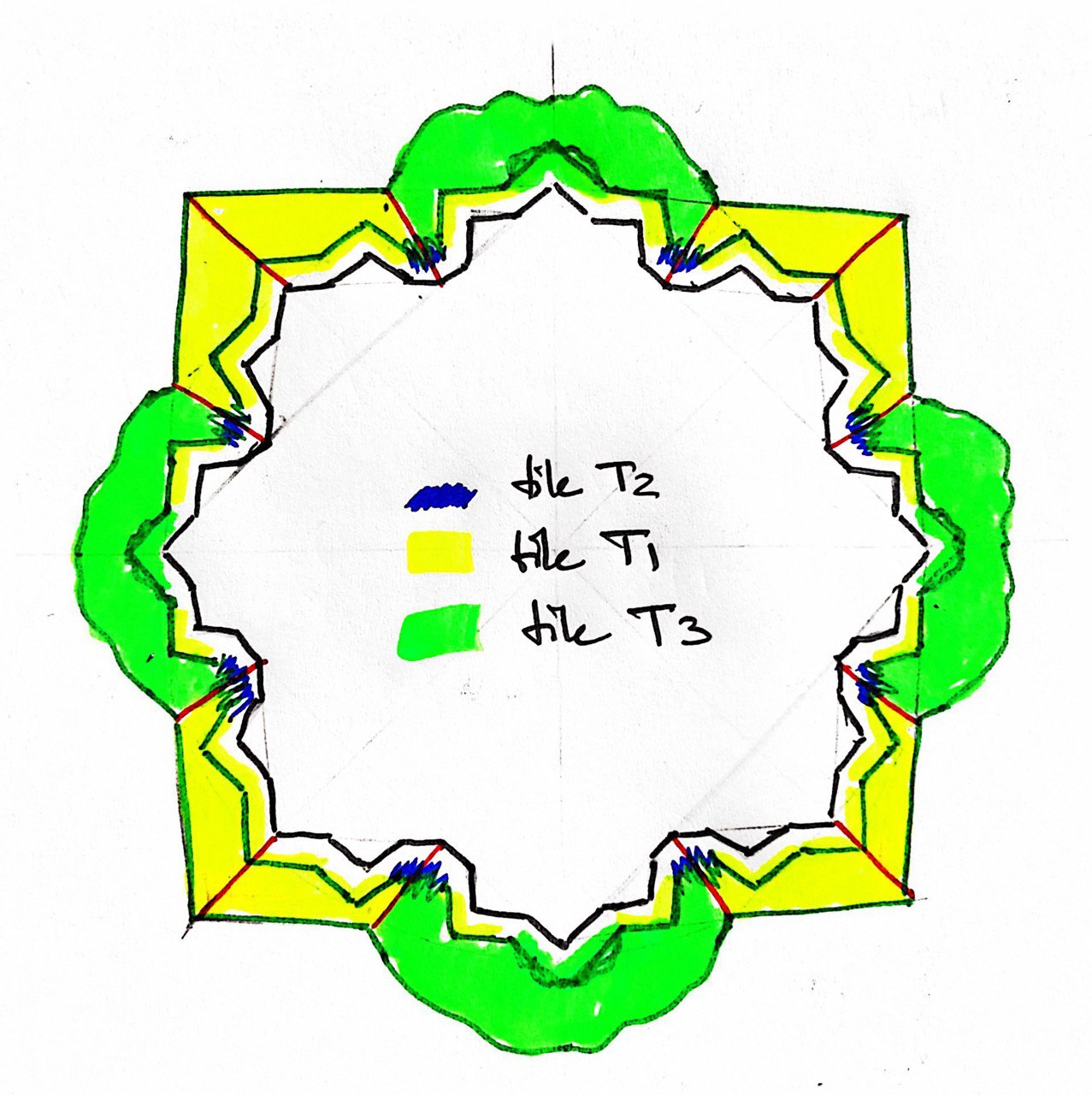}
    \caption{Tiling for a compact snowflake with any $\alpha$}
    \label{fig:comp tiling any}
\end{figure}

We can do versions without the restriction $\alpha = \frac{\pi}{N}, N \ge 4 \in \N$, as well, but this is less elegant. For example, we could start with a regular polygon and treat all its sides as initial segments $I = [z_1, z_2]$, build snowflakes on them, tile their neighbourhoods like we did in Section 4, and then do some gluing with tiles of another type around vertices of the initial polygon like we did for our non-compact versions of the snowflakes in the previous section. Once again, we just give a sketch f tiling around the snowflake, see figure $\ref{fig:comp tiling any}$.


\begin{figure}[h]
    \centering
    \includegraphics[scale = 0.2]{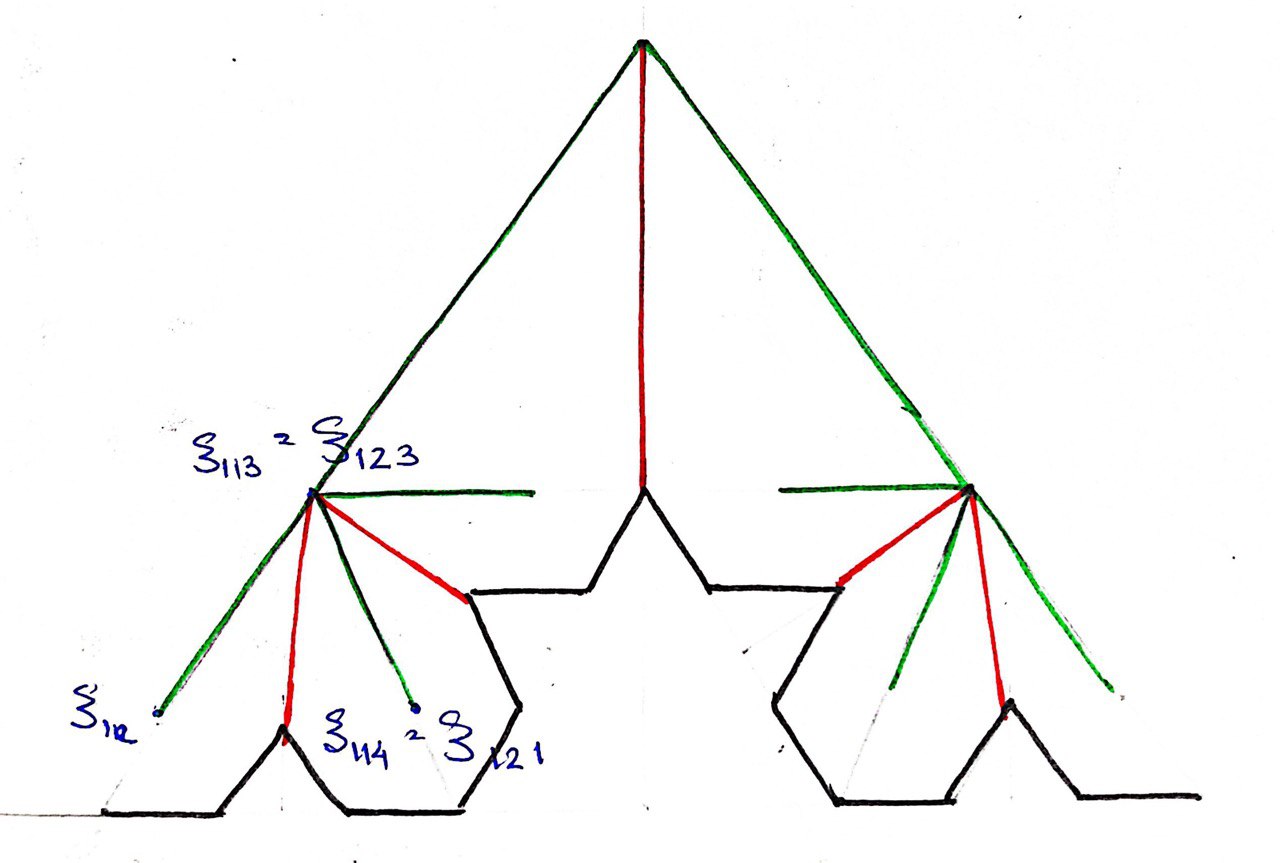}
    \caption{The limit case of $\alpha = \pi/3$}
    \label{fig:lim case}
\end{figure}

The case of the classical Koch snowflake with the angle $\alpha = \pi/3$ is not covered by Theorem $\ref{mainthm1}$. ``Going to a limit'' using precisely the construction we have above does not work for the following reason. Consider the curve $s_2$ which we use later as a model to construct all other approximating curves $s_k$. As $\alpha$ goes to $\pi/3$, points $\zeta_{113}$ and $\zeta_{123}$, as well as $\zeta_{114}$ and $\zeta_{122}$, get closer and closer to each other. This means that at $\alpha = \pi/3$ segments $[\zeta_{113}, \zeta_{114}]$ and $[\zeta_{122}, \zeta_{123}]$ are the same, so different parts of our curve are glued together, and instead of a proper part of a curve we have three green segments $[\zeta_{112}, \zeta_{113}], [\zeta_{113}, \zeta_{114}]$ and $[\zeta_{123}, \zeta_{124}]$ which stick out from the same point $\zeta_{113} = \zeta_{123}$. See figure $\ref{fig:lim case}$.

This also implies that the green sides of our tile $T2$ get more and more distorted, and that red sides of tiles get shorter and shorter: in fact, their length tends to zero when $\alpha$ tends to $\pi/3$. Since the constant $C$ in the inequality $\eqref{scalar_coeff}$ we aim for can be roughly estimated as the ratio of the length of an upper green lid of a tile and the length of a red lid of a tile, we see that the constant $C = C_\alpha$ indeed tends to infinity as $\alpha$ goes to $\pi/3$, as mentioned in the introduction. 

We suspect, however, that, by twitching our construction slightly, one could cover this case as well. While defining $\eqref{s1}$, we are not really obliged to choose this particular form of the curve $(\zeta_{12}, \zeta_{14})$. We could do something more flat, avoiding the gluing described above for the limiting case. 

~\

Concerning the open questions we have, it would be interesting to see if our construction could be adapted for the case when, instead of building an isosceles triangle always above the segment we are modifying, one replaces the middle of the segment of two sides of isosceles triangle build either above or below the segment, depending on where we are in the construction. The same question but about varying the angle $\alpha$ of the isosceles triangles depending on where we are (probably combined with the sign change) is also amusing. We suspect that one cannot do something completely random. For example, changing randomly the sign of the angle $\alpha$ in the construction, which corresponds to the choice of whether to build a triangle above of below the base segment, should not be allowed. That is, one can probably still show that (actually we suspect this really can be done, see later), for small enough absolute values of $\alpha$, if the sign of the angle changes randomly, there still exists a good operator in the sense of Theorem $\ref{mainthm1}$. We just doubt that such a result can be obtained by the methods this paper describes. To get a good operator with our tiling method, some regularity in the structure of a step of the construction of a snowflake should be maintained, so that, at the end of the day, we are able to have a finite set of tiles which covers the whole upper half-space (or at least a band adjacent to its boundary).

We wonder if our method is applicable to other self-similar fractals in $\R^2$. There are plenty of fractal sets which are generated by actions of some mappings like our favourite transformation $F$ above. The positive answer could widen the range of dimensions for which a good operator exists, as lots of those sets have dimensions larger than $\frac{\ln{(4)}}{\ln{(3)}}$. Here we would like to share only an impression that comes from some drawing experiments. It seems like, the higher dimension of the fractal is, the harder it is to pack nice approximating curves in the space (because the fractal itself has to be packed rather densely to have a large dimension). 

Returning to the three questions we announced in the introduction, we think that one can significantly extend the class of examples of unrectifiable sets for which a good operator in sense of Theorem $\ref{mainthm1}$ exists. Namely, we conjecture that for all Reifenberg flat sets on the plane with a small enough (flatness) constant one can construct such an operator. For the definition and some examples of Reifenberg flat sets, see, for example, \cite{DT1}. Note that these sets can be very irregular, so we doubt that some tiling technique can work in this case. One would have to invent a completely new approach.

Finally, we would like to build some examples of good operator for ``bad'' sets in higher dimensions. We can do some already, as mentioned in Section 5 of \cite{DM}, by taking sets like $S_\alpha \times \R$ in $\R^3$ and an accompanying operator with the coefficient $A(x, y, z) = a(x, y)$, but this does not look like a proper higher dimensional fractal. The main difficulty here is to invent the procedure which reconstructs an operator coefficient from the level surfaces of its solution. To our knowledge, such a procedure, presented in Section 2, is known only in dimension 2. It also seems that identifying the correct set of restrictions on the coefficient $A$ so that they are geometrically more relevant to the story than just a coefficient obtained from a generalization of a quasiconformal mapping from a half-space above the higher dimensional snowflake to the standard half-space (in dimension 2 we have a scalar-valued coefficient versus a matrix-valued one), might be not easy. And there once again we wonder if one can construct a good operator for all higher-dimensional Reifenberg flat sets with a small enough constant, as opposed to some snowflake-type sets with a rich group of symmetries.

~\

\noindent Polina Perstneva\\
Universit\'e Paris Saclay, LMO \\
polina.perstneva@universite-paris-saclay.fr

\end{document}